\documentclass[10pt]{amsart}

\usepackage[colorlinks]{hyperref}
\usepackage{color,graphicx,shortvrb}
\usepackage[latin 1]{inputenc}

\usepackage[active]{srcltx} 

\usepackage{enumerate}
\usepackage{amssymb}

\usepackage{amsmath}
\usepackage{amsfonts}
\usepackage{amsthm}
\usepackage{mathtools}
\usepackage{xcolor}

\newtheorem{theorem}{Theorem}[section]

\newtheorem{corollary}[theorem]{Corollary}

\newtheorem{lemma}[theorem]{Lemma}

\newtheorem{proposition}[theorem]{Proposition}
\newtheorem{remark}[theorem]{Remark}

\newcommand*{\myproofname}{Proof of the Claim}
\newenvironment{inproof}[1][\myproofname]{\begin{proof}[#1]}{\end{proof}}

\def\J#1#2#3{ \left\{ #1,#2,#3 \right\} }
\def\RR{{\mathbb{R}}}
\def\NN{{\mathbb{N}}}
\def\11{\textbf{$1$}}
\def\CC{{\mathbb{C}}}

\newcommand{\R}{\mathbb{R}}
\newcommand{\N}{\mathbb{N}}

\newcommand{\eps}{\varepsilon}
\newcommand{\dd}{{\ast\ast}}
\newcommand{\spn}{\mathrm{span}}
\newcommand{\df}{\coloneqq}

\newcommand{\zuz}{\!\!\upharpoonright}
\newcommand{\AAA}{{\mathfrak{A}}}
\newcommand{\B}{{\mathcal{B}}}

\newcommand{\jed}{\mathbf{1}}

\DeclareGraphicsExtensions{.jpg,.pdf,.png,.eps}

\begin{document}

\title[Kadison's transitivity theorem for unital JB$^*$-algebras and the MUP]{A strengthened Kadison's transitivity theorem for unital JB$^*$-algebras with applications to the Mazur--Ulam property}

\author[A.M. Peralta]{Antonio M. Peralta}

\address{Instituto de Matem{\'a}ticas de la Universidad de Granada (IMAG). Departamento de An{\'a}lisis Matem{\'a}tico, Facultad de Ciencias, Universidad de Granada, 18071 Granada, Spain.}
\email{aperalta@ugr.es}

\author[R. \v{S}varc]{Radovan \v{S}varc}

\address{Department of Mathematical Analysis, Faculty of Mathematics and Physics, Charles University,	Sokolovská 83, 186 00 Prague, Czech Republic}
\email{Rada.svarc@seznam.cz}




\date{December 24th, 2022}

\begin{abstract} The principal result in this note is a strengthened version of Kadison's transitivity theorem for unital JB$^*$-algebras, showing that for each minimal tripotent $e$ in the bidual, $\AAA^{**}$, of a unital JB$^*$-algebra $\AAA$, there exists a self-adjoint element $h$ in $\AAA$ satisfying $e\leq \exp(ih)$, that is, $e$ is bounded by a unitary in the principal connected component of the unitary elements in $\AAA$. This new result opens the way to attack new geometric results, for example, a Russo--Dye type theorem for maximal norm closed proper faces of the closed unit ball of $\AAA$ asserting that each such face $F$ of $\AAA$ coincides with the norm closed convex hull of the unitaries of $\AAA$ which lie in $F$. Another geometric property derived from our results proves that every surjective isometry from the unit sphere of a unital JB$^*$-algebra $\AAA$ onto the unit sphere of any other Banach space is affine on every maximal proper face. As a final application we show that every unital JB$^*$-algebra $\AAA$ satisfies the Mazur--Ulam property, that is, every surjective isometry from the unit sphere of $\AAA$ onto the unit sphere of any other Banach space $Y$ admits an extension to a surjective real linear isometry from $\AAA$ onto $Y$. 
This extends a result of M. Mori and N. Ozawa who have proved the same for unital C$^*$-algebras. 
\end{abstract}

\maketitle
\thispagestyle{empty}

\section{Introduction}

\emph{Kadison's transitivity theorem} is a starring result in the theory of C$^*$-algebras which surprisingly asserts the equivalence of the notions of topological irreducibility and algebraic irreducibility for representations of C$^*$-algebras. More concretely,  let $A$ be a non-zero C$^*$-algebra acting
irreducibly on a Hilbert space $H$, and suppose that $\xi_1, \ldots, \xi_n$ and $\eta_l, \ldots, \eta_n$
are elements in $H$ such that $\xi_1, \ldots, \xi_n$ are linearly independent. Then there
exists an operator $a \in A$ such that $a(\xi_j) = \eta_j$ for $j = 1,\ldots, n$. If there is
a self-adjoint operator $b$ in $B(H)$ such that $b(\xi_j) = \eta_j$ for $j = 1,\ldots, n$, then
we may choose $a$ to be self-adjoint in $A$. If A contains the identity operator on $H$ and there is a
unitary $u$ in $B(H)$ such that $a(\xi_j) = \eta_j$ for $j = 1,\ldots, n$, then we may choose $a$
to be a unitary -- we may even suppose that $a = \exp(iw)$ for some self-adjoint element $w \in A_{sa}$ (cf. \cite[Theorem 5.4.3]{KadRingrBook1}, \cite[Theorem 3.13.2]{Ped} or \cite[Theorem 5.2.2]{MurphyBook}). An equivalent reformulation affirms that if $p_1,\ldots, p_n$ are mutually orthogonal minimal projections in the bidual, $A^{**}$, of a C$^*$-algebra $A$, then there exist mutually orthogonal norm-one positive elements $a_1,\ldots,a_k$ in $A$ such that $a_j = p_j + (\mathbf{1}-p_j) a_j (\mathbf{1}-p_j)$ for all $j = 1,\ldots, n$ (note that $(\mathbf{1}-p_j) a_j (\mathbf{1}-p_j)$ is value of the Peirce-$0$ projection of $p_j$ applied to $a$). The result was extended to the setting of JB$^*$-algebras by J. Hamhalter in \cite[[Proposition 2.3]{Ham99}.\smallskip
 
A more non-associative approach leads us to a version of Kadison's transitivity theorem for JB$^*$-triples, due to L.J. Bunce, F.J. Fern{\'a}ndez-Polo, J. Mart{\'i}nez and the first author of this note, which asserts that for every finite family of mutually orthogonal minimal tripotents $e_1,\ldots, e_n$ in the bidual, $E^{**}$, of a JB$^*$-triple $E$, we can find mutually orthogonal norm-one elements  $a_1, \ldots , a_n$ in $E$ such that $a_j = e_j + P_0(e_j) (a_j)$ for $j = 1, \ldots , n$ \cite[Theorem 3.3]{BuFerMarPe}. The detailed definitions of JB$^*$-algebras and JB$^*$-triples, as well as the basic terminology, are all gathered in subsection \ref{subsec:background}. \smallskip

Kadison's transitivity theorem, and specially the statement concerning unital C$^*$-algebras and unitary elements, is a crucial tool in a recent result by M. Mori and N. Ozawa (cf. \cite[\S 5]{MoriOza2018}), where they prove that every unital C$^*$-algebra $A$ satisfies the Mazur--Ulam property, that is  every surjective isometry from the unit sphere of $A$ onto the unit sphere of any other Banach space extends to a surjective real linear isometry between the spaces. The key refinement shows that for each minimal partial isometry $e$ in $A^{**}$ there exists a unitary element $u$ in $A$ such that $e \leq u$ with respect to the natural order among partial isometries (i.e. $u-e$ is partial isometry and $u-e$ is orthogonal to $e$). One of the main goals in this note, and probably the most striking, is to establish a domination result of this type in the wider setting of unital JB$^*$-algebras.\smallskip

In support of our research proposal we remark a recent contribution on finite JB$^*$-algebras and JB$^*$-triples in \cite{HKP-Fin}; a notion of finiteness which coincides with the usual concept for projections in von Neumann algebras. According to the just quoted reference, every minimal tripotent in a JBW$^*$-triple (i.e. a JB$^*$-triple which is also a dual Banach space) is finite \cite[Lemma 3.2$(e)$]{HKP-Fin}, and for each finite tripotent $e$ in a JBW$^*$-algebra $\mathfrak{M}$ there is a unitary element $u\in \mathfrak{M}$ with $e\leq u$ with respect to the natural order on tripotents \cite[Proposition 7.5]{HKP-Fin}.\smallskip

We establish here a strengthened version of Kadison's transitivity theorem for unital JB$^*$-algebras showing that for each minimal tripotent $e$ in the bidual, $\AAA^{**}$, of a unital JB$^*$-algebra $\AAA$, there exists a self-adjoint element $h$ in $\AAA$ satisfying $e\leq \exp(ih)$, that is, $e$ is bounded by a unitary in the principal connected component of the unitary elements in $\AAA$ (see Theorem~\ref{t boundedness  principle for minimal tripotents in the second dual by unitaries in the principal component}). The proof of this result follows a completely new approach which is based on a detailed study of the JB$^*$-subalgebra generated by a minimal tripotent in a general JB$^*$-algebra (cf. Proposition~\ref{p representation of the JBstar algebra generated by a minimal projection}), and a subtle combination of the classical Kadison's transitivity theorem for JB$^*$-triples \cite{BuFerMarPe} and the recent results on finite JBW$^*$-algebras \cite{HKP-Fin}. \smallskip

The lacking of polar decompositions for elements in JBW$^*$-algebras has been a handicap to extend classical result on von Neumann algebras and C$^*$-algebras to the Jordan setting. Our strengthened version of Kadison's transitivity theorem for unital JB$^*$-algebras opens a new way to attack new geometric results without any use of polar decompositions. For example, since the facial structure of the closed unit ball, $\mathcal{B}_{E},$ of a general JB$^*$-triple $E$ relies on the knowledge of compact tripotents in its second dual \cite{EdFerHosPe2010}, and particularly, the maximal norm closed proper faces of $\mathcal{B}_{E}$ are in one-to-one correspondence with minimal tripotents in $E^{**},$ we can find a first application of the result on a Russo--Dye type theorem for maximal norm closed proper faces of the closed unit ball of a unital JB$^*$-algebra $\AAA$. Concretely, each maximal norm closed proper face $F$ of $\AAA$ coincides with the norm closed convex hull of the unitaries of $\AAA$ which lie in $F$ (see Theorems~\ref{Russo_Dye_theorem_for_F_p} and \ref{Russo_Dye_theorem_for_maximal proper faces}), which is a improved version of the Russo--Dye theorem for unital JB$^*$-algebras \cite{WriYou78,Sidd2010,Sidd2006,Sidd2007,MackMell21}.\smallskip

The paper contains some other geometric tools designed, like a study of the convex combination of a proper norm closed face of the closed unit ball of a JB$^*$-triple with its antipodal (see section~\ref{sec: convex combination of faces with their antipodals}), or a technical results proving that every surjective isometry from the unit sphere of a unital JB$^*$-algebra $\AAA$ onto the unit sphere of any other Banach space is affine on every maximal proper face (see Proposition~\ref{Prop_20_from_Mori-Ozawa} and Theorem~\ref{theorem prop_20_from_Mori-Ozawa for maximal faces associated with minimal tripotents}). Although a Russo--Dye type theorem for general real JBW$^*$-algebras remains as an open problem, we prove that every convex body $\mathcal{K}$ in a real JBW$^*$-algebra satisfies the strong Mankiewicz property in the sense defined in \cite{MoriOza2018},, that is, every surjective isometry from $\mathcal{K}$ onto an arbitrary convex subset $L$ in a normed space $Y$ is affine (cf. Proposition~\ref{p strong Mankiewicz for the real JB$^*$-algebras associated with a Jordan multiplicative functional}, where an more general conclusion is obtained).\smallskip

We have been naturally led to one of the most striking open problems in functional analysis, \emph{Tingley's problem} on the extension of surjective isometries between the unit spheres of two normed spaces to surjective real linear isometries between the whole spaces (see, for example, \cite{Ting1987, Ding2009, YangZhao2014, FangWang06, ChenDong2011, Pe2018, Tan2014, Tan2016, Tan2017b, PeTan16, JVMorPeRa2017, Pe2019, FerPe17c, FerPe17b, FerPe17d, FerPe18Adv, FerGarPeVill17, FerJorPer2018, CuePer19, CabSan19,  Mori2017, MoriOza2018, BeCuFerPe2018, KalPe2019,  Hat2021, Banakh21, BanakhCabello21, Pe2022Tingleycompact, Hat2023} to realise the intense research activity around this problem). Tingley's problem remains open even for finite dimensional Banach spaces of dimension greater than or equal to three, a positive solution when one of the involved spaces is 2-dimensional has been found by T. Banakh in 2022 \cite{Banakh22}. Let us recall that a Banach space $X$ satisfies the \emph{Mazur--Ulam property} if Tingley's problem admits a positive answer when $X$ is one of the involved spaces \cite{ChenDong2011}. M. Mori and N. Ozawa established in \cite{MoriOza2018} that every unital C$^*$-algebra satisfies the Mazur--Ulam property. The main results in \cite{BeCuFerPe2018,KalPe2019} prove that the same conclusion holds for JBW$^*$-triples (i.e. JB$^*$-triples which are dual Banach spaces), in particular for JBW$^*$-algebras. The question whether the Mazur--Ulam property is automatically satisfied by any unital JB$^*$-algebra (without assuming duality) remains as a natural problem until now. The study of this problem has required some years to develop the geometric tools we are presenting in this paper. The final section of this manuscript is devoted to show that every unital JB$^*$-algebra satisfies the Mazur--Ulam property (see Theorem \ref{t unital JB$^*$-algebras satisfy the MUP}).

\subsection{Background on JB$^*$-triples, the facial structure of the closed unit ball, and the strong Mankiewicz property }\label{subsec:background}\ \smallskip

JB$^*$-algebras constitute a deeply studied mathematical model, which provides a non-associative extension of C$^*$-algebras, and shares most of the geometric properties of the latter. The pioneering contribution by P. Jordan, J. von Neumann and E. Wigner \cite{JorvNeuWign34} introduced Jordan algebras with the aim of being a suitable setting for axiomatic 
quantum mechanics. The notions of JB- and JB$^*$-algebras began to be studied from the point of view of functional analysis during the mid-1960s and 1970s (see the monographs \cite{HOS,AlfsenShultz2003,CabGarPalVol1}). A \emph{Jordan algebra} over the field $\mathbb{K}$ ($\mathbb{R}$ or $\mathbb{C}$) is a non-necessarily associative algebra $\AAA$ whose (Jordan) product, denoted by $\circ$, satisfies $x\circ y = y \circ x$ and the so-called \emph{Jordan identity}: $$( x \circ y ) \circ x^2 = x\circ ( y\circ x^2 ) \hbox{ for all } x,y\in \AAA.$$ Jordan algebras are \emph{power 
	associative}, that is, each subalgebra generated by a single element $a$ is associative, equivalently, by setting $a^0 = \textbf{1}$, and $a^{n+1} = 
a\circ a^n,$ we have $a^n \circ a^m = a^{n+m}$ for all $n,m\in \mathbb{N}\cup \{0\}$. (cf. \cite[Lemma 2.4.5]{HOS} or \cite[Corollary 1.4]{AlfsenShultz2003}). Given an element $a \in \mathfrak{A}$ the symbol $U_a$ will stand for the linear mapping on $\mathfrak{A}$ defined by $U_a (b) := 2(a\circ b)\circ a - a^2\circ b.$ 
\smallskip

A Jordan-Banach algebra $\mathfrak{A}$ is a Jordan algebra equipped with a complete norm satisfying $\|a\circ b\|\leq \|a\| \cdot \|b\|$ for all $a,b\in \mathfrak{A}$. A JB$^*$-algebra is a complex Jordan-Banach algebra $\mathfrak{A}$ equipped with an algebra involution $^*$ satisfying the following Jordan version of the Gelfand-Naimark axiom: $$\|U_a ({a^*}) \|= \|a\|^3, \ (a\in \mathfrak{A}).$$ 
Let us observe that every C$^*$-algebra is a JB$^*$-algebra with respect to the natural Jordan product $x\circ y = \frac12 ( xy + yx)$, the same involution and norm. In this particular case $U_a ({a^*}) = aa^* a,$ and hence $\| aa^* a\| =\|U_a ({a^*}) \|= \|a\|^3$ is an equivalent reformulation of the Gelfand--Naimark axiom. The JB$^*$-subalgebras of C$^*$-algebras are known as \emph{JC$^*$-algebras}. The class of JC$^*$-algebras is strictly smaller than the class of JB$^*$-algebras; there exist exceptional JB$^*$-algebras which cannot be embedded as Jordan $^*$-subalgebras of a C$^*$-algebra (see, for example,  \cite[Corollary 2.8.5]{HOS}).\smallskip
 
A \emph{JB-algebra} is a real Jordan algebra $\mathfrak{J}$ equipped with a complete norm satisfying \begin{equation}\label{eq axioms of JB-algebras} \|a^{2}\|=\|a\|^{2}, \hbox{ and } \|a^{2}\|\leq \|a^{2}+b^{2}\|\ \hbox{ for all } a,b\in \mathfrak{J}.
\end{equation} 

The indissoluble connection between JB- and JB$^*$-algebras was already posed by I. Kaplansky during his famous lecture to the 1976 St. Andrews Colloquium of the Edinburg Mathematical Society. The key contribution by J.D.M. Wright showed that every JB-algebra $\mathfrak{J}$ corresponds uniquely to the self-adjoint part $\mathfrak{A}_{sa}=\{x\in \mathfrak{A} : x^* =x\}$ of a JB$^*$-algebra $\mathfrak{A}$ \cite{Wri77}. A JBW$^*$-algebra (respectively, a JBW-algebra) is a JB$^*$-algebra (respectively, a JB-algebra) which is also a dual Banach space. Every JBW$^*$-algebra (respectively, each JBW-algebra) contains a unit element (see \cite[\S 4]{HOS} or \cite{AlfsenShultz2003}). It is worth to note that JBW-algebras are precisely the self-adjoint parts of JBW$^*$-algebras (see \cite[Theorems 3.2 and 3.4]{Ed80JBW} or \cite[Corollary 2.12]{MarPe00}).\smallskip 

Several basic notions in JB$^*$-algebra theory are better understood from the wider perspective of JB$^*$-triples. As shown in \cite{Ka83}, there are certain biholomorphic properties (namely, being a bounded symmetric domain) of the open unit ball of a complex Banach space which are only satisfied by the open unit ball of those complex Banach spaces in the class of JB$^*$-triples. A \emph{JB$^*$-triple} is a complex Banach space $E$ admitting a continuous triple product $\J \cdot\cdot\cdot : E\times E\times E \to E,$ which is symmetric and bilinear in the outer variables, conjugate-linear in the inner one, and satisfies the following axioms:
\begin{enumerate}[{\rm (a)}] \item (Jordan identity)
	$$L(a,b) L(x,y) = L(x,y) L(a,b) + L(L(a,b)x,y)
	- L(x,L(b,a)y)$$ for $a,b,x,y$ in $E$, where $L(a,b)$ is the operator on $E$ given by $x \mapsto \J abx;$
	\item $L(a,a)$ is a hermitian operator with non-negative spectrum for all $a\in E$;
	\item $\|\{a,a,a\}\| = \|a\|^3$ for each $a\in E$.\end{enumerate}

For each element $a\in E$ the symbol $Q(a)$ will denote the conjugate linear operator on $E$ defined by $Q(a) (x) = \{a,x,a\}$ ($x\in E$).\smallskip

A JBW$^*$-triple is a JB$^*$-triple which is also a dual Banach space. In such a case, it admits a unique isometric predual and its triple product is separately weak$^*$-continuous (cf. \cite{BarTi86}).\smallskip 

Harris \cite[Corollary 2]{Harris74} proved that the open unit ball of a C$^*$-algebra $A$ is a bounded symmetric domain. Actually, $A$ is a JB$^*$-triple for the triple product given by \begin{equation}\label{eq triple product JCstar triple} \{a,b,c\} = \frac12 (a b^* c + c b^* a) \ \ (a,b,c\in A).
\end{equation} The triple product just commented is also valid to induce a structure of JB$^*$-triple on the space $B(H,K),$ of all bounded linear operators between two complex Hilbert spaces $H$ and $K$. JB$^*$-triples of the form  $B(H,K)$ are called \emph{Cartan factor of type 1}. The remaining five types of Cartan factors are described next. Cartan factors of types 2 and 3 are subtriples of $B(H)$ similar to symmetric and skew-symmetric complex matrices. Fix a conjugation $j$ (i.e. a conjugate-linear isometry or period 2) on a complex Hilbert space $H$, and define a linear involution on $B(H)$ by $x\mapsto x^t:=jx^*j$ --this is an infinite dimensional version of the transposition in $M_n(\mathbb{C})$. \emph{Cartan factors of type 2 and 3} are the JB$^*$-subtriples of $B(H)$ of all $t$-skew-symmetric and $t$-symmetric operators, respectively.\smallskip

A \emph{Cartan factor of type 4}, also known as a \emph{spin factor},\label{def spin factor} is a complex Hilbert space $M$ provided with a conjugation $x\mapsto \overline{x},$ where the triple product and the norm are defined by \begin{equation}\label{eq spin product}
	\{x, y, z\} = \langle x, y\rangle z + \langle z, y\rangle  x -\langle x, \overline{z}\rangle \overline{y},
\end{equation} and \begin{equation}\label{eq spin norm} \|x\|^2 = \langle x, x\rangle  + \sqrt{\langle x, x\rangle ^2 -|
		\langle x, \overline{x}\rangle  |^2},
\end{equation} respectively (cf. \cite[Chapter 3]{Fri2005}). The \emph{Cartan factors of types 5 and 6} (also called \emph{exceptional} Cartan factors) are spaces of matrices over the eight dimensional complex algebra of Cayley numbers; the type 6 consists of all $3\times 3$ self-adjoint matrices and has a natural structure of JB$^*$-algebra, and the type 5 is the subtriple consisting of all $1\times 2$ matrices which is not a JB$^*$-algebra (see \cite{Ka97, Harris74} and the recent references \cite[\S 6.3 and 6.4]{HKP-Fin}, \cite[\S 3]{HamKalPe22determinants} for more details).\smallskip  

A central result in JB$^*$-triple theory, known as \emph{Kaup's Banach--Stone theorem} asserts that a linear bijection between JB$^*$-triples is an isometry if and only if it is a triple isomorphism (cf. \cite[Proposition 5.5]{Ka83}).\smallskip

From the point of view of the arguments in this paper, it is worth noticing that every JB$^*$-algebra $\AAA$ is a JB$^*$-triple with respect to the triple product defined by $$\{a,b,c\} := (a\circ b^*)\circ c + (c\circ b^*)\circ a - (a\circ c) \circ b^* \ \  \ (a,b,c\in \AAA).$$ As we shall see in section~\ref{sec: Russo--Dye for faces in unital JBstaralgebras}, when we deal with unitaries in unital JB$^*$-algebras, each unitary element in a unital JB$^*$-algebra $\AAA$ induces a Jordan product and an involution defining a new structure of unital JB$^*$-algebra on $\AAA$. Although different unitaries in $\AAA$ can produce JB$^*$-algebra structures which are not Jordan $^*$-isomorphic (cf. $(\mathcal{PU}2)$ in page~\pageref{PU2}), by Kaup's Banach--Stone theorem, there is only one triple product on the Banach space $\AAA$ (see $(\mathcal{PU}3)$ in page~\pageref{PU3}).\smallskip  

Back to the triple product employed to regard each C$^*$-algebra $A$ as a JB$^*$-triple (cf. \eqref{eq triple product JCstar triple}), we realise that an element $e\in A$ is a partial isometry if and only if $\{e,e,e\} =e$. An element $e$ in a general JB$^*$-triple $E$ is called a \emph{tripotent} if $\{e,e,e\} =e$.  Each tripotent $e$ in $E$ induces a decomposition (known as the \emph{Peirce decomposition} of $E$ associated with $e$) of the space $E$ in terms of the eigenspaces of the operator $L(e,e)$:
\begin{equation}\label{Peirce decomp} {E} = {E}_{0} (e) \oplus  {E}_{1} (e) \oplus {E}_{2} (e),\end{equation} where ${
	E}_{k} (e) := \{ x\in {E} : L(e,e)x = {\frac k 2} x\}$ is a subtriple of ${E}$ called the \emph{Peirce-$k$ subspace} ($k=0,1,2$). \emph{Peirce-$k$ projection} is the name given to the natural projection of ${E}$ onto ${E}_{k} (e)$ and it is usually denoted by $P_{k} (e)$.
Triple products among elements in Peirce subspaces obey the following patterns: the inclusion $\J {{E}_{k}(e)}{{E}_{l}(e)}{{E}_{m}(e)}\! \subseteq {E}_{k-l+m} (e),$ and the identity $\J {{E}_{0}(e)}{{E}_{2} (e)}{{E}}\! =\! \J {{E}_{2} (e)}{{E}_{0} (e)}{{E}}\! =\! \{0\},$ hold for all $k,l,m\in \{0,1,2\}$, where ${E}_{k-l+m} (e) = \{0\}$ whenever $k-l+m$ is not in $\{0,1,2\}$. We shall employ the term \emph{Peirce arithmetic} to refer to these rules.\smallskip  

It is perhaps worth to recall that the projection of $E$ onto $E_j(e)$ (denoted by $P_j(e)$ and called the Peirce-$j$ projection) is a contractive mapping for all $j=0,1,2$ (cf. \cite{FriRu85}). Peirce projections obey the following identities $P_{2}(e) = Q(e)^2,$ $P_{1}(e) =2(L(e,e)-Q(e)^2),$ and $P_{0}(e) =Id_E - 2 L(e,e) + Q(e)^2.$\smallskip

The Peirce-$2$ subspace ${E}_{2} (e)$ is a unital JB$^*$-algebra with respect to the product and involution given by $x \circ_e y = \J xey$ and $x^{*_e} = \J exe,$ respectively. 
\smallskip

A non-zero tripotent $e$ in a JB$^*$-triple $E$ is called (\emph{algebraically}) \emph{minimal} if  $E_2(e)=\CC e \neq \{0\}$. The tripotents $e\in E$ satisfying $E_0(e) =\{0\}$ (respectively, $E_2(e) = E$) are called \emph{complete} (respectively, \emph{unitaries}). When a C$^*$-algebra $A$ is regarded as a JB$^*$-triple, minimal partial isometries, maximal partial isometries and unitaries correspond to minimal, complete and unitary tripotents, respectively.\smallskip

It is known that the extreme points of the closed unit ball of a JB$^*$-triple $E$ are precisely the complete tripotents in $E$, that is \begin{equation}\label{eq extreme points and complete tripotents} \partial_e(\mathcal{B}_{E}) =\{ \hbox{complete tripotents in } E \},
\end{equation} (cf. \cite[Lemma 4.1]{BraKaUp78} and \cite[Proposition 3.5]{KaUp77}).\smallskip

A (closed) subtriple $I$ of a JB$^*$-triple $E$ is said
called an \emph{ideal} or a \emph{triple ideal} of $E$ if $\J EEI + \J EIE\subseteq I$. If we can only guarantee that the containing $\J IEI \subseteq I$ holds, we say that $I$ is an \emph{inner ideal} of $E$. \smallskip

JB$^*$-triples admit a Gelfand theory for the subtriple generated by a single element $a$ in a JB$^*$-triple $E$. If we write $a^{[1]} := a$, $a^{[3]} := \{a,a,a\}$, and $a^{[2n+1]} := \J a{a^{[2n-1]}}a,$ $(n\in \NN)$, and we denote by $E_a$ the JB$^*$-subtriple of $E$ generated by $a$, it is known that there exists an isometric triple isomorphism $\Psi_a$ from $E_a$ onto $C_0 (\Omega_{a})$ for some compact Hausdorff space $\Omega_{a}$ contained in $[0,\|a\|],$ such that $\|a\|\in \Omega_{a}$, where $C_0 (\Omega_{a})$ denotes the Banach space of all complex-valued continuous functions vanishing at $0$ if $0\in \Omega_a$, such that $\Psi_a (a)$ corresponds to the embedding of $\Omega_a$ into $\mathbb{C}$ (cf. \cite[Lemma 1.14]{Ka83} or \cite{Ka96}). The set $\Omega_a$ is known as the \emph{triple spectrum} of $a$. Thanks to this local Gelfand theory, for each $f\in C_0 (\Omega_{a})$, we write $f_t(a) = \Psi_a^{-1} (f)$. We define this way the triple functional calculus of $f$ at the element $a$. As an application, given a natural $n$, there exists (a unique) $a^{[{1}/({2n-1})]}$ in
$E_a$ satisfying $(a^{[{1}/({2n-1})]})^{[2n-1]} = a.$ For any $a\in E$ the sequence $(a^{[{1}/({2n-1})]})$ converges in the weak$^*$ topology of $E^{**}$ to a tripotent, denoted by $r(a),$ called the \emph{range tripotent} of $a$ in $E^{**}$. The tripotent $r(a)$ is the smallest tripotent $e$ in $E^{**}$ satisfying that $a$ is positive in the JBW$^*$-algebra $E^{**}_{2} (e)$. It is also known that the sequence $(a^{[2n -1]})$ converges in the weak$^*$ topology of $E^{**}$ to a tripotent (called the \hyphenation{support}\emph{support} \emph{tripotent} of $a$) $u(a)$ in $E^{**}$, which satisfies $ u(a) \leq a \leq r(a)$ in $E^{**}_2 (r(a))$ (compare \cite[Lemma 3.3]{EdRu88}).\smallskip

Orthogonality is another concept required in our arguments. Recall that elements $a,b$ in a JB$^*$-triple $E$ are said to be \emph{orthogonal} (written $a\perp b$) if $L(a,b) =0$. It is known \cite[Lemma 1]{BurFerGarMarPe} that $a\perp b$ if and only if any of the following statements holds:
\begin{equation}\label{eq reformulations of orthogonality}
\begin{array}{ccc}
	a\perp b; & \J aab =0; & a \perp r(b);  \\
	& & \\
	r(a) \perp r(b); & E^{**}_2(r(a)) \perp E^{**}_2(r(b)); & r(a) \in E^{**}_0 (r(b)); \\
	& & \\
	a \in E^{**}_0 (r(b)); & b \in E^{**}_0 (r(a)); & E_a \perp E_b. \\
\end{array}
\end{equation}	

A subset $S$ of a JB$^*$-triple $E$ is called orthogonal if it does not contain zero and satisfies that $x \perp y$ for every $x\neq y$ in $S$. The \emph{rank of $E$} is the minimal cardinal number $r$ satisfying $\hbox{card}(S) \leq r$ for every \emph{orthogonal} subset $S \subseteq E$. The reader can consult the references \cite{Ka97}, \cite{BunChu92} and \cite{BeLoPeRo} for the basic theory on the rank of a Cartan factor and of a JB$^*$-triple, and its relation with reflexivity.\smallskip

\noindent The Peirce rules assure that for each tripotent
$u$ in a JB$^*$-triple $E$, the sets $E_0(u)$ and $E_2 (u)$ are orthogonal.\smallskip

There is a natural order among tripotents in a JB$^*$-triple $E$ defined in terms of the relation of orthogonality which is given by $e\leq v$ if and only if $v-e$ is a tripotent orthogonal to $e$ (cf. \cite{HKP-Fin,HamKalPe22orders}).\smallskip 

Let us consider a weak$^*$ dense subset $E$ in the dual $X^*$ of a Banach space $X$. A non-zero subset $\mathcal{O}$ of $X^{*}$ is called \emph{open relative to $E$} if $\mathcal{O}\cap E$ is $\sigma(X^{*},X)$-dense in $\overline{\mathcal{O}}^{\sigma(X^{*},X)}$ (c.f.
\cite{FerPe06}). Let $E$ be a weak$^*$ dense JB$^*$-subtriple of a
JBW$^*$-triple $W$. A tripotent $u$ in $W$ is said to be
\emph{open} relative to $E$ if $W_2 (u)$ is open relative to $E$ (compare \cite[p. 167]{EdRu96} and \cite[p. 78]{FerPe06}). A tripotent $u$ in $W$ such that $W_0 (u)$ is open relative to $E$ is called \emph{closed} relative to $E$. Finally, a tripotent $u$ is called \emph{bounded} relative to $E$ if there exists a norm-one element $x$ in $E$ such that $x = u + P_0 (u) (x)$ (c.f. \cite{FerPe06}).\smallskip

The notion of compact tripotent in $W$ relative to $E$ is due to Edwards and Ruttimann \cite[\S 4]{EdRu96}. Following \cite{EdRu96}, we shall say that a tripotent $u$ in $W$ is \emph{compact} relative to $E$ if $u=0$ or there exists a decreasing net, $(u_{\lambda}) \subseteq W,$ of support tripotents associated to norm-one elements in $E$, converging in the weak$^*$ topology of $W,$ to $u$. By \cite[Theorem 2.6]{FerPe06} a tripotent $u\in W$ is compact relative to $E$ if and only if $u$ is closed and bounded relative to $E$.\smallskip

The notion of compactness for tripotents is consistent with the one existing for projections and partial isometries in the C$^*$-algebra setting. Furthermore, if a C$^*$-algebra $A$ is regarded as a JB$^*$-triple, then a
projection (respectively, a partial isometry) $p$ in $A^{**}$ is open relative to $A$ if and only if $p$ is an open tripotent relative to $A$ (compare \cite{Ak69} or \cite{AkPed1,AkPed92} or \cite[Proposition 3.11.9]{Ped}).\smallskip

Compact tripotents are the natural notion to understand the facial structure of the closed unit ball of a JB$^*$-triple. Let $F$ and $G$ be two subsets of the closed unit ball, $\mathcal{B}_{X},$ of a Banach space $X$ and the closed unit ball, $\mathcal{B}_{X^*}$, of its dual, respectively. We define 
$$ F^{\prime} =F^{\prime,X^*} =  \{a \in \mathcal{B}_{X^*}:a(x) = 1\,\, \forall x \in F\},$$
$$G_{\prime} =G_{\prime,X}= \{x \in \mathcal{B}_{X}:a(x) = 1\,\, \forall a \in G\}.$$
Clearly, $F^{\prime}$ is a weak$^*$ closed face of $\mathcal{B}_{X^*}$ and $G_{\prime}$
is a norm-closed face of $\mathcal{B}_{X}$. We say that $F$ is a \emph{norm-semi-exposed face} of $\mathcal{B}_{X}$ (respectively, $G$ is a \emph{weak$^*$ semi-exposed face} of $\mathcal{B}_{X^*}$) if $F=(F^{\prime})_{\prime}$ (respectively, $G= (G_{\prime})^{\prime}$). It is known that the mappings $F \mapsto F^{\prime}$
and $G \mapsto G_{\prime}$ are anti-order isomorphisms between the
complete lattices $\mathcal{S}_n(\mathcal{B}_{X})$, of norm-semi-exposed faces
of $\mathcal{B}_X,$ and $\mathcal{S}_{w^*}( \mathcal{B}_{X^*}),$ of weak$^*$ semi-exposed
faces of $\mathcal{B}_{X^*}$, and are inverses of each other.\smallskip

If $M$ is a JBW$^*$-triple, Edwards and Ruttimann proved in \cite{EdRu88} that every weak$^*$ closed
face of $\mathcal{B}_{M}$ is weak$^*$ semi-exposed; furthermore, the mapping \begin{equation}\label{eq anti order weak* closed faces complex} u \mapsto
	(\{u\}_{\prime})^{\prime} = u + \mathcal{B}_{M_0(u)}
\end{equation} is an anti-order isomorphism from the set of all tripotents in $M$ onto the complete lattice of all weak$^*$ closed faces of $\mathcal{B}_{M}$ (see \cite[Theorem 4.6]{EdRu88}). Another result in the just quoted reference proves that if $M_*$ stands for the predual of $M$, every norm-closed face of $\mathcal{B}_{M_{*}}$ is norm-semi-exposed, and the mapping
\begin{equation}\label{eq order isomorphism norm-closed faces predual complex} u \mapsto \{u\}_{\prime}
\end{equation} is an order isomorphism from the set of all tripotents in $M$ onto the complete lattice of all norm-closed faces of $\mathcal{B}_{M_{*}}$ \cite[Theorem 4.4]{EdRu88}.\smallskip

In 1992, Akemann and Pedersen described all norm-closed faces of the closed unit ball of a C$^*$-algebra $A$ and the weak$^*$ closed faces of $\mathcal{B}_{A^*}$ in terms of compact partial isometries in $A^{**}$ (see \cite{AkPed92}). In 2010, Edwards, Fernández-Polo, Hoskin and the first author of this note established a description of the norm-closed faces of the closed unit ball of a JB$^*$-triple -- the result reads as follows:

\begin{theorem}\label{t facil structure of the ball of a JB$^*$-triple}{\rm(\cite[Corollaries 3.11 and 3.12]{EdFerHosPe2010} and \cite[Corollary 3.5 and the preceding comments]{BuFerMarPe})} Let $E$ be a JB$^*$-triple. Then every norm-closed face of $\mathcal{B}_{E}$ is norm-semi-exposed and the mapping 
	$$u \mapsto F_{u} :=(\{u\}_{\prime})_{\prime} = (u + \mathcal{B}_{E^{**}_0(u)})\cap E$$
	is an anti-order isomorphism from the set of all compact tripotents in $E^{**}$ onto the complete lattice of all norm-closed faces of $\mathcal{B}_{E}$. Furthermore, each maximal proper face of $\mathcal{B}_{E}$ is of the form $F_e$ for a suitable minimal tripotent in $E^{**}$.
\end{theorem}

An additional geometric aspect of JBW$^*$-triples is required for our purposes. The extreme points, $\partial_e\left(\mathcal{B}_{W_*}\right)$, of the closed unit ball, $\mathcal{B}_{W_*}$, of the predual, $W_*$, of a JBW$^*$-triple $W$ are called pure atoms\label{pure atoms and support minimal tripotents} (cf. \cite{FriRu85}).  For each pure atom $\varphi\in  \partial_e\left(\mathcal{B}_{W_*}\right)$, there exists a unique minimal tripotent $e=e_{\varphi}\in W$ (called the \emph{support minimal tripotent} of $\varphi$) satisfying $P_2(e) (x) = \varphi(x) e$ for all $x\in W$, and reciprocally, the Peirce-2 projection associated with each minimal tripotent in $W$ is given by an extreme point of the closed unit ball of $W_*$ (see \cite[Proposition 4]{FriRu85}). A minimal tripotent $e\in W$ is the support tripotent of a pure atom $\varphi \in \partial_e\left(\mathcal{B}_{W_*}\right)$ if and only if $\varphi (e)=1$.  It is also known from \cite{BurFerGarMarPe,EdFerHosPe2010} that for each JB$^*$-triple $E$, each maximal proper face of of $\mathcal{B}_{E}$ is of the form $\{\varphi\}_{\prime}$ for a unique $\varphi\in  \partial_e\left(\mathcal{B}_{E^*}\right)$ and $\{\varphi\}_{\prime} = (e_{\varphi} + \mathcal{B}_{E^{**}_0(e_{\varphi})})\cap E$ where $e_{\varphi}$ denotes the support minimal tripotent of $\varphi$ in $E^{**}$. In case that $E = \AAA$ is a JB$^*$-algebra, the extreme points of the positive part of $\mathcal{B}_{\AAA^*}$ are known as pure states of $\AAA$ and are in one-to-one correspondence with the minimal projections in $\AAA^{**}$.\smallskip

For the purposes of this note, we recall that L. Chen and Y. Dong were the first authors observing that each surjective isometry $\Delta$ between the unit spheres of two Banach spaces maps maximal convex subsets to maximal convex subsets \cite[Lemma 5.1]{ChenDong2011}, equivalently, $\Delta$ preserves maximal proper faces of the corresponding closed unit balls (cf. \cite[Lemma 3.3]{Tan2016}, \cite[Lemma 3.5]{Tan2014}). \smallskip

In the case of a JB$^*$-triple $E$, every proper norm closed face of $\mathcal{B}_E$ is norm-semi-exposed (cf. Theorem \ref{t facil structure of the ball of a JB$^*$-triple} or \cite[Corollary 3.11]{EdFerHosPe2010}). Among the arguments in the proof of \cite[Proposition 2.4]{FerGarPeVill17} we can find a proof of the fact that every norm-semi-exposed face of $\mathcal{B}_E$ coincides with the intersection of all maximal proper norm closed faces containing it, that is, every proper norm closed face of $\mathcal{B}_E$ is an intersection face according to the notation in \cite{MoriOza2018}. Actually Lemma 8 in \cite{MoriOza2018} proves that every surjective isometry $\Delta$ between the unit spheres of two Banach spaces maps intersection faces to intersection faces. We gather all these conclusion in the next lemma.

\begin{lemma}\label{l intersection faces MoriOzawa L8 and semi-exposed faces FerGarPeVill}{\rm(\cite[Lemma 8]{MoriOza2018}, \cite[Proposition 2.4 and Corollary 2.5]{FerGarPeVill17})} Let $\Delta: S_E\to S_Y$ be a surjective isometry where $E$ is a JB$^*$-triple and $Y$ is a real Banach space. Then $\Delta$ maps proper norm closed faces of $\mathcal{B}_{E}$ to intersection faces in $S_Y$. Furthermore, if $F$ is a proper norm closed face of $\mathcal{B}_{E}$ then $\Delta(-F) = -\Delta(F)$.
\end{lemma}

Proposition 2.4 in \cite{FerGarPeVill17} actually proves that every surjective isometry between the unit spheres of two Banach spaces satisfying an additional geometric property preserves norm-semi-exposed faces in the sphere. 

\section{Convex combinations of faces of the closed unit ball of a JB$^*$-triple with their antipodals}\label{sec: convex combination of faces with their antipodals}

This section is devoted to understand the geometric--algebraic connections which allow us to characterize the convex combination of a proper norm closed face of the closed unit ball of a JB$^*$-triple and its antipodal. Henceforth all faces will be assumed to be closed. Let $F$ be a proper face of the closed  unit ball of a JB$^*$-triple $E$. Given $\lambda\in[-1,1]$ we set $$ F_{\lambda} := \Big\{ x\in S_{E} : \hbox{dist}(x,F)\leq 1-\lambda\ \& \  \hbox{dist}(x,-F)\leq 1+\lambda  \Big\}.$$ This ``geometric'' device is employed successfully applied in \cite{MoriOza2018} (also in \cite{CuePer19,Hat2021,HatOiTog21,CabCuHiMiPe22,Hat2023}), where it is observed that, since $\hbox{dist}(F,-F) = 2$, the inequalities $\leq$ in the previous definition can be replaced with equalities. In the very particular case in which $p$ is a compact projection in the second dual of a C$^*$-algebra $A$,  and $F = F_p =A\cap \left( p + \mathcal{B}_{(\mathbf{1}-p) A^{**} (\mathbf{1}-p)}\right)$ is the proper face of $\mathcal{B}_{A}$ associated with $p$, Mori and Ozawa proved in \cite[Lemma 17]{MoriOza2018} that $$\left(F_p\right)_{\lambda} = F(p,\lambda) := S_{A}\cap \left( \lambda p + \mathcal{B}_{(1-p) A^{**} (1-p)}\right)$$ for all $\lambda\in [-1,1]$. The question whether this algebraic--geometric identity holds for all compact partial isometries in $A^{**}$, or more generally, for all compact tripotents in the second dual of a general JB$^*$-triple $E$ remains open. This is our main task in this section.

\begin{theorem}\label{t lemma 17 for triples} Let $e$ be a compact tripotent in the bidual of a JB$^*$-triple $E$, and let $F_e = E\cap (e+\mathcal{B}_{E_0^{**}(e)})$ denote the proper face of $\mathcal{B}_E$ associated with $e$. Then for each $\lambda\in [-1,1]$ we have $$  \left(F_e\right)_{\lambda} = F(e,\lambda) = S_{E}\cap \left( \lambda e + \mathcal{B}_{E^{**}_0 (e)}\right).$$
\end{theorem}

The proof will be obtained after a series of lemmata. We begin with some arguments inspired from \cite[Lemma 17]{MoriOza2018}.

\begin{lemma}\label{l distance to pm1 Cstar} Let $a$ be an element in a unital C$^*$-algebra $A$ such that $$\| a - \mathbf{1}\| + \| a+ \mathbf{1} \| =2.$$ Then there exits a real number $\lambda$ in $[-1,1]$ such that $a = \lambda \mathbf{1}$.
\end{lemma} 

\begin{proof} We can assume, via Gelfand--Naimark theorem, that $A \subset B(H)$ for some complex Hilbert space $H$. Given any $\xi$ in the unit sphere of $H$ we have 
$$ 2 = 2\|\xi\| \leq \| (\textbf{1} - a) (\xi) \| + \| (a+ \textbf{1}) (\xi)\| \leq \| a - \textbf{1}\| + \| a+ \textbf{1} \| =2.$$ It then follows that there exists a real number $\lambda_{\xi}\in [-1,1]$, depending on $\xi$, such that  $a(\xi) = \lambda_{\xi} \xi$. If $H$ is one-dimensional there is nothing to prove, otherwise, take any other element $\eta$ in the unit sphere of $H$.\smallskip

First suppose that $\eta\not\perp\xi$. Then we can find $\alpha,\beta\in \mathbb{C}$ with $\alpha\ne0$ and $\eta'$ in the unit sphere of $H$ orthogonal to $\xi$ such that $\eta = \alpha \xi + \beta \eta'$ and $|\alpha|^2 + |\beta|^2 =1$.  It follows from what is proved in the first part that there exist $\lambda_{\xi},\lambda_{\eta},\lambda_{\eta'}\in [-1,1]$ satisfying the previously commented identity, and hence $$ \lambda_{\eta} \alpha  \xi + \lambda_{\eta} \beta  \eta' = \lambda_{\eta} \eta = a(\eta) = \alpha a(\xi) + \beta a(\eta') =  \alpha \lambda_{\xi} \xi + \beta \lambda_{\eta'} \eta'.$$ Since $\alpha\ne0$, this shows that $\lambda_{\eta} = \lambda_{\xi}$.\smallskip

Finally, if $\eta\perp\xi$, we can denote $\zeta\df\frac{\xi+\eta}{\|\xi+\eta\|}$. Since $\eta\not\perp\zeta\not\perp\xi$, we get $\lambda_\eta=\lambda_\zeta=\lambda_\xi$ from the previous part. Therefore $\lambda_\xi$ doesn't depend on $\xi$.
\end{proof}

Let $a$ be any element in a unital JB$^*$-algebra $\mathfrak{A}$. We recall that the celebrated Shirshov-Cohn theorem \cite[Theorem 2.4.14]{HOS} assures that the JB$^*$-subalgebra of $\mathfrak{A}$ generated by $a$ and $\mathbf{1}$ is isometrically JB$^*$-isomorphic to a JB$^*$-subalgebra subalgebra of a unital C$^*$-algebra $A$, and we can even more assume that $\mathbf{1}$ is the unit of $A$. The next corollary is thus a consequence of the previous conclusion and these comments.

\begin{corollary}\label{c distance to pm1 JBstar} Let $a$ be an element in a unital JB$^*$-algebra $\mathfrak{A}$ such that $$\| a - \mathbf{1}\| + \| a+ \mathbf{1} \| =2.$$ Then there exits a real number $\lambda$ in $[-1,1]$ such that $a = \lambda \mathbf{1}$.
\end{corollary}  

Let us recall some properties of Peirce subspaces associated with a tripotent $e$ in a JB$^*$-triple $E$. As we have already commented, the Peirce-2 subspace $E_2(e)$ is a unital JB$^*$-algebra with unit $e,$ Jordan product $a\circ_e b = \{a,e,b\}$ and involution $a^{*_e}= \{e,a,e\}$. It follows from Peirce arithmetic that $$P_2 (e) \{x,x,e\} = \{P_2(e)(x),P_2(e)(x),e\} + \{P_1(e)(x),P_1(e)(x),e\},$$ for all $x\in E$. It is further known that $\{P_2(e)(x),P_2(e)(x),e\}$ and  $\{P_1(e)(x),$ $P_1(e)(x),e\}$ are positive elements in $E_2(e)$ and $\{P_j(e)(x),P_j(e)(x),e\} =0$ for some $j=1,2$ if and only if $P_j(e) (x) =0$ (\cite[Lemma 1.5]{FriRu85} and \cite{Pe15}). Let $\mathcal{S} (E_2(e))$ denote the set of all states on $E_2(e)$ (i.e. the set of all positive norm-one functionals in $E_2(e)^*$). Since the norm of every positive element in a unital JB$^*$-algebra $\mathfrak{A}$ is attained at its evaluation at a state on $\mathfrak{A}$, the assignment  $$ x\mapsto \|x\|_e := \sup\left\{ \phi\{x,x,e\}^{\frac12} = (\phi P_2(e) \{x,x,e\})^{\frac12} : \phi \in \mathcal{S} (E_2(e))\right\}$$ defined a Hilbertian norm on $E_2(e)\oplus E_1(e)$ which is equivalent to the original one (cf. \cite[Proposition 2.4, Corollary 2.5 and comment prior to the latter]{BurFerGarMarPe}). If $E$ is a JBW$^*$-triple, we can replace the set $\mathcal{S} (E_2(e))$ of all states on $E_2(e)$ by the set $\mathcal{S}_n (E_2(e))$ of all normal states on the JBW$^*$-algebra $E_2(e)$ and everything remains invariant. Let us be more explicit, the results in the just commented reference show that for $x= x_2 +x_1\in E_2(e)\oplus E_1(e)$ we have \begin{equation}\label{eq first inequlity towards equivalence of norms} 
	\|x\|_e\leq \|x\|\leq \|x_2\| + \|x_1\| \leq \|x_2\|_e + 2 \|x_1\|_e .
\end{equation} The projections $P_j(e)$ ($j=0,1,2$) are norm contractive for the original norm  (cf. \cite[Corollary 1.2]{FriRu85}). We can actually see that for $j=1,2$, the restriction $P_j (e)|_{E_2(e) \oplus E_1(e)}: E_2(e) \oplus E_1(e) \to E_j(e)$ is contractive for the norm $\|\cdot\|_e$ for all $j=1,2$. Indeed, by the comments above we have 
$$\begin{aligned}
	 \|P_j (e) (P_2 (e) (x) &+ P_1 (e) (x) )\|^2_e =\sup_{\phi \in \mathcal{S} (E_2(e)) } \phi\{P_j (e) (x) , P_j (e) (x)  ,e\} \\ &\leq \sup_{\phi \in \mathcal{S} (E_2(e)) } \phi\{P_2 (e) (x) , P_2 (e) (x)  ,e\} +\phi\{ P_1 (e) (x) , P_1 (e) (x)  ,e\} \\
	&=\sup_{\phi \in \mathcal{S} (E_2(e)) } \phi\{P_2 (e) (x) + P_1 (e) (x) , P_2 (e) (x) + P_1 (e) (x)  ,e\} \\ &=  \| P_2 (e) (x) + P_1 (e) (x) \|^2_e,
\end{aligned} $$  for all $j=1,2,$ $x \in E$. By combining the previous conclusion with \eqref{eq first inequlity towards equivalence of norms} we derive that 
\begin{equation}\label{eq equivalence of norms} 
	\|x\|_e\leq \|x\|\leq 3 \|x\|_e, \hbox{ for all } x\in E_2(e)\oplus E_1(e).
\end{equation}

Given a functional $\varphi$ in the dual of a JB$^*$-triple $E$ and a tripotent $e\in E$ with $\|\varphi\| = \varphi (e)$ we know that $\varphi = \varphi P_2(e)$ \cite[Proposition 1]{FriRu85}. \smallskip

The Hilbertian norm $\|\cdot\|_e$ is also related to the preHilbertian seminorm in the Grothendieck's inequalities for JB$^*$-triples (cf. \cite{BarFri87Groth,HKPP-BF}). Given a normal functional $\varphi$ in the predual of a JBW$^*$-triple $W$, and any $z\in S_W$ with $\|\varphi\| = \varphi(z)$, the mapping $x\mapsto \|x\|_{\varphi}:=\varphi\{x,x,z\}^{\frac12}$ ($x\in W$) is a preHilbertian seminorm on $W$, which does not depend on the chosen $z\in S_W$. If $\phi$ is any normal state on a JBW$^*$-algebra $\mathfrak{A}$, we have $\|x\|_{\phi}^{2} = \phi\{x,x,\mathbf{1}\} = \phi (x\circ x^*)$ for all $x\in \mathfrak{A}$. According to this notation, for each tripotent $e$ in a JB$^*$-triple $E$ we have $$\|x\|_e = \sup\left\{ \|x\|_{\phi} : \phi \in \mathcal{S} (E_2(e))\right\},$$ where we identify the states on $E_2(e)$ with the normal states on $E_2^{**}(e)$. It is also known (cf. \cite[\S 3]{BarFri90} or \cite[Proposition 1.2]{BarFri87Groth}) that given a norm-one functional $\varphi$ in the dual of a JB$^*$-triple $E$, we have \begin{equation}\label{eq bound of norm varphi} |\varphi (x) |\leq \|x\|_{\varphi}, \hbox{ for all } x\in E.
\end{equation} 

We can now understand the true meaning of Lemma~\ref{l distance to pm1 Cstar} and Corollary~\ref{c distance to pm1 JBstar} in the case of tripotents in general JB$^*$-triples. 

\begin{lemma}\label{l distance to pm tripotent JBtriple} Let $a$ be an element in a JB$^*$-triple $E$, and let $e$ be a tripotent in $E$ such that $\| a - e\| + \| a+ e \| =2$. Then there exits a real number $\lambda\in[-1,1]$ such that $a = \lambda e + P_0(e)(a)$.
\end{lemma}

\begin{proof} Let us begin with the inequality 
	$$\begin{aligned}2 = \| 2 e\| \leq \| e - P_2(e) (a) \| + \| P_2(e) (a) +e\| &= \| P_2(e) ( e - a) \| + \| P_2(e) (a +e) \|\\ 
		&\leq \| e - a\| + \| a+ e \| =2,
	\end{aligned}$$  which proves that $2 = \| e - P_2(e) (a) \| + \| P_2(e) (a) +e\|$, and thus Corollary~\ref{c distance to pm1 JBstar} applied to the element $P_2(e)(a)$ in the unital JB$^*$-algebra $E_2(e)$, gives the existence of a real number $\lambda\in [-1,1]$ such that $P_2(e) (a) = \lambda e$. \smallskip

On the other hand, the mapping $L(e,e)$ is contractive and due to $L(e,e) = P_2(e) +\frac 12 P_1(e)$ satisfies $P_2(e)=P_2(e)L(e,e)$. Therefore  $$\begin{aligned}2 = \| 2 e\| \leq \| e - P_2(e) (a) \| + \| P_2(e) (a) +e\| &\leq  \| L(e,e) ( e - a) \| + \| L(e,e) (a +e) \|\\ 
	&\leq \| e - a\| + \| a+ e \| =2,
\end{aligned}$$ which guarantees that $$ \| e - P_2(e) (a) \| =  \| L(e,e) ( e - a) \| = \| L(e,e) ( a) -e \|$$ and $$ \| e + P_2(e) (a) \| =  \| L(e,e) ( e + a) \| = \| L(e,e) ( a) +e \|.$$ By applying that $P_2(e) (a) = \lambda e$ for some $\lambda\in [-1,1]$, we deduce that $$ |1-\lambda| = \| e - \lambda e \| =  \left\| (1-\lambda) e - \frac12 P_1(e) (a) \right\| $$ and $$ |1+\lambda| = \| e + \lambda e \| =  \left\| (1+\lambda) e + \frac12 P_1(e) (a) \right\|.$$ If $\lambda = \pm 1$, the above identities give $P_1(e) (a) =0$. We can therefore assume that $-1 < \lambda <1$, and hence $\displaystyle 1 =  \left\| e - \frac{1}{2(1-\lambda) } P_1(e) (a) \right\| $. We fix now an arbitrary state $\phi\in \mathcal{S} (E_2(e))$ and check the following inequality $$\begin{aligned} 1 & = \phi \left( e - \frac{1}{2(1-\lambda) } P_1(e) (a) \right)^2 \leq \left\| e - \frac{1}{2(1-\lambda) } P_1(e) (a) \right\|_{\phi} ^2\\ &= \left\| e  \right\|_{\phi} ^2 + \left\|  \frac{1}{2(1-\lambda) } P_1(e) (a) \right\|_{\phi} ^2
 = 1 + \left\|  \frac{1}{2(1-\lambda) } P_1(e) (a) \right\|_{\phi} ^2 \\ &=  \left\| e - \frac{1}{2(1-\lambda) } P_1(e) (a) \right\|_{\phi} ^2 \leq \left\| e - \frac{1}{2(1-\lambda) } P_1(e) (a) \right\|^2 =1\end{aligned} ,$$ where we applied \eqref{eq bound of norm varphi} and Peirce arithmetic. The previous inequalities show that $$ 0 = \left\|  \frac{1}{2(1-\lambda) } P_1(e) (a) \right\|_{\phi} ^2 = \phi \left\{\frac{1}{2(1-\lambda) } P_1(e) (a),\frac{1}{2(1-\lambda) } P_1(e) (a),e \right\},$$ for all states $\phi \in   \mathcal{S} (E_2(e))$, which combined with the positivity of the element $\left\{\frac{1}{2(1-\lambda) } P_1(e) (a),\frac{1}{2(1-\lambda) } P_1(e) (a),e \right\}$ in $E_2(e)$, gives $$\left\{\frac{1}{2(1-\lambda) } P_1(e) (a),\frac{1}{2(1-\lambda) } P_1(e) (a),e \right\} =0,$$ and thus $P_1(e) (a) =0$ by \cite[Lemma 1.5]{FriRu85} or \cite{Pe15}.
\end{proof}

One ingredient is missing to prove the main result of this section is a generalized Urysohn's lemma taken from \cite{FerPeUrysohn}.

\begin{proof}[Proof of Theorem~\ref{t lemma 17 for triples}] Let us first take $x\in  \left(F_e\right)_{\lambda}$, then for each $\varepsilon>0$ there exist $y= e+ P_0(e) (y)\in F_e$ and $z= -e+ P_0(e) (z)\in -F_e$ satisfying $\|x- y \| <1-\lambda +\frac{\varepsilon}{2}$ and $\|x- z \| <1+\lambda+\frac{\varepsilon}{2} $.  Having in mind that $L(e,e) = P_2(e)+\frac12 P_1(e)$ is a contractive mapping on $E^{**}$ we derive that 
$$\begin{aligned} 2 = \| 2 e \| &\leq \| e - L(e,e) (x) \| + \| e + L(e,e) (x) \| \\
	& = \| L(e,e) (y - x) \| + \| L(e,e) (z - x) \| \\ 
	&\leq  \|y-x\| + \| z-x\| < 2 + \varepsilon. 
\end{aligned} $$ It follows from the arbitrariness of $\varepsilon>0$ that $2 = \| e - L(e,e) (x) \| + \| e + L(e,e) (x) \| $. Lemma~\ref{l distance to pm tripotent JBtriple} applied to the tripotent $e$ and the element $L(e,e) (x) $ in $E^{**}$ assures the existence of a real number $\lambda \in [-1,1]$ satisfying $L(e,e) (x) = \lambda e$. This implies that $x = \lambda e + P_0(e) (x) \in F(e,\lambda) = S_{E}\cap \left(\lambda e + \mathcal{B}_{E_0^{**}(e)} \right)$ --the Peirce decomposition is taken in $E^{**}$.\smallskip

Take now $x \in F(e,\lambda) = S_{E}\cap \left(\lambda e + \mathcal{B}_{E_0^{**}(e)} \right).$ If $\lambda =0$, we have $x = P_0 (e) (x)$ is orthogonal to $e$. Fix $\varepsilon>0.$ Let us consider the function $f:[0,1]\to \mathbb{R}$ given by
$$f (s):=\left\{ \begin{array}{ll}
	0 , & \hbox{if $0\leq s\leq \varepsilon$} \\
	\hbox{affine }, & \hbox{ if $\varepsilon \leq s\leq 2 \varepsilon$} \\
	s, & \hbox{ if $2\varepsilon \leq s\leq 1$}. 
\end{array} \right.$$ the element $f_t (x)\in E$ given by the continuous triple functional calculus, and the compact tripotent $w = \chi_{[\varepsilon, 1]} (x) \in E^{**}$. By the non-commutative Urysohn's lemma for orthogonal compact tripotents established in \cite[Corollary 3.6 or Proposition 3.7]{FerPeUrysohn}, there exists a norm-one elements $a\in E$ such that $a = e+P_0(e) (a)$ and $a\perp w$.  Clearly, since $f_t(x)\in E^{**}_2(w)$, the element $z = a+ f_t(x)$ lies in $F_e,$ and by construction $$\hbox{dist} (x,F_e) \leq \| z - x\| \leq \| a \| + \|f_t(x) -x\| = 1 + \varepsilon.$$ The arbitrariness of $\varepsilon>0$ proves that $\hbox{dist} (x,F_e) \leq 1 = 1-\lambda$. Similarly, $y = - a+ f_t(x)$ lies in $F_{-e},$ and by construction $$\hbox{dist} (x,F_{-e}) \leq \| y - x\| \leq \| - a \| + \|f_t(x) -x\| = 1 + \varepsilon,$$ witnessing that  $\hbox{dist} (x,-F_e) \leq 1 = 1+\lambda$.\smallskip

Suppose now that $\lambda \neq 0$. We can assume that $0<\lambda \leq 1$, otherwise we replace $e$ with $-e$. Take another function $g:[0,1]\to \mathbb{R}$ defined by
$$g (s):=\left\{ \begin{array}{ll}
	0 , & \hbox{if $s = 0$} \\
	\hbox{affine }, & \hbox{ otherwise,} \\
	1, & \hbox{ if $s\in[\lambda,1]$}. 
\end{array} \right.$$ The element $g_t(x)\in S_{E}$. Since $x = \lambda e + P_0(e) (x)$, by orthogonality, it is easy to check that $g_t (x) = g_t(\lambda e) + g_t (P_0(e)(x)) = e + g_t(P_0(e)(x))\in F_e$ (computed in $E^{**}$), and $\| x - g_t(x)\| = 1-\lambda$. This proves that $\hbox{dist}(x,F_e) \leq 1-\lambda$. We can similarly get $\hbox{dist}(x, -F_e) = \hbox{dist}(x, F_{-e}) \leq 1+\lambda$.\smallskip

Alternatively, if $x = \lambda e + P_0 (e) (x) \in F(e,\lambda),$ the elements $y = e + P_0(e) (x)$ and $z= - e + P_0(e) (x) $ lie in $F_e^{E^{**}}= e + \mathcal{B}_{E_0(e)}$ and $-F_e^{E^{**}} =F_{-e}^{E^{**}} = -e + \mathcal{B}_{E_0(e)}$, respectively. Clearly, $\|x-y\|= \|(1-\lambda)e\|=  1-\lambda$ and $\|x-z\| =  1+\lambda$, and hence $$\hbox{dist}\left(x, F_e^{E^{**}}\right)\leq 1-\lambda\ \& \  \hbox{dist}\left(x,-F_e^{E^{**}}\right)\leq 1+\lambda.$$  By \cite[Theorem 3.6]{BeCuFerPe2018}, $F_e^{E^{**}}$ (respectively, $-F_e^{E^{**}}$) coincides with the weak$^*$ closure of the proper face $F_e = S_E\cap F_e^{E^{**}}$ (respectively, $-F_e = S_E\cap -F_e^{E^{**}}$), that is, $F_e^{E^{**}} = \overline{F_e}^{\sigma(E^{**},E^*)}$ (respectively, $F_{-e}^{E^{**}} = \overline{F_{-e}}^{\sigma(E^{**},E^*)}$). We claim that, by the Hahn-Banach separation theorem  $ \hbox{dist}(x, F_e) = \hbox{dist}\left(x, F_e^{E^{**}}\right)$ and $ \hbox{dist}(x, -F_e) = \hbox{dist}\left(x,- F_e^{E^{**}}\right)$. Namely, if $ \hbox{dist}(x, F_e) =0$ there is nothing to prove. Otherwise, there exists a functional $\varphi\in S_{E^*}$ such that $\displaystyle \sup_{a\in F_e} \Re\hbox{e} \varphi(a) +\hbox{dist} (x, F_e) \leq   \Re\hbox{e} \varphi(x)$ (cf. \cite[2.2.19 Mazur's Separation Theorem or 2.2.26 Eidelheit's Separation Theorem]{Megg98}). Since $\varphi$ lies in the predual of $E^{**},$  we have $\sup_{a\in F_e} \Re\hbox{e} \varphi(a) = \sup_{z\in \overline{F_e}^{\sigma(E^{**},E^*)}} \Re\hbox{e} \varphi(z)$, and hence $$ \hbox{dist} (x, F_e) \leq \| x-z\|, \ \hbox{ for all } z\in \overline{F_e}^{\sigma(E^{**},E^*)},$$ which proves the claim. 
\end{proof}

In subsequent sections we shall also need the following technical result borrowed from \cite{MoriOza2018}

\begin{lemma}\label{l 11 MoriOzawa}{\rm\cite[Lemmata 9 and 11]{MoriOza2018}} Let $\Delta: S_{X}\to  S_{Y}$ be a surjective isometry between the unit spheres of two Banach spaces. Then the following statements hold:\begin{enumerate}[$(a)$] 
		\item $\Delta$ maps intersection faces in $S_{X}$ to intersection faces in $S_{Y}$. 
		\item The identity $\Delta(F_{\lambda}) = \Delta(F)_{\lambda}$ holds for every  intersection face $F\subseteq S_{X}$ and every $\lambda\in [-1,1]$. 
		\item Let $F$ be an intersection face in $S_{X}.$  Then, for any $\lambda_1,\lambda_2$ in $[-1,1]$ and $\alpha\in [0, 1]$ we have 
		$$(\alpha\Delta(F_{\lambda_1}) + (1-\alpha) \Delta(F_{\lambda_2}))\cap S_{Y} \subseteq \Delta \left( F_{\alpha \gamma_1 + (1-\alpha)\gamma_2}\right).$$
	\end{enumerate}

\end{lemma}

\section{Russo--Dye theorem for faces}\label{sec: Russo--Dye for faces in unital JBstaralgebras}

The celebrated Russo--Dye theorem asserts that the closed unit ball of a unital C$^*$-algebra $A$ coincides with the norm-closed convex hull of the unitaries in $A$ \cite{RuDye, Gard84}. Actually every element in the open unit ball of $A$ is the mean of a finite number of unitaries in $A$ \cite{KadPed}. L. Harris introduced in \cite{Harris72} holomorphic theory to establish a variant of the Russo--Dye theorem for unital complex Banach algebras with involution, denoted by $^*$, where the notion of unitary element is the natural one. It was proved that the closed convex hull of the unitary elements of these algebras contains all elements $x$ whose spectral radius and the spectral radius of $x^* x$ is smaller than or equal to $1$, if the inequalities are strict then the element belongs to the convex hull of the unitaries. There is another strengthened version of this result, proved by M. Mori and N. Ozawa, asserting that every maximal proper face $F$ of the closed unit ball of $A$ coincides with the closed convex hull of all unitary elements in $F,$ that is, $F = \overline{co} \left( F\cap \mathcal{U} (A)\right)$ \cite[Lemma 18]{MoriOza2018}. \smallskip

One of the many uses of the Russo--Dye theorem became more tangible after the following theorem due to Mori and Ozawa. We recall first that a convex subset in a normed space $X$ is called a
\emph{convex body} if it has non-empty interior in $X$.

\begin{theorem}\label{t Mori-Ozawa strong Mankiewicz property}{\rm\cite[Theorem 2]{MoriOza2018}}. Let $X$ be a Banach space such that the closed convex hull of the extreme points $\partial_e \left(\mathcal{B}_X\right)$, of the closed unit ball, $\mathcal{B}_X$, of $X$ has non-empty interior in $X$. Then every convex body $\mathcal{K} \subseteq  X$ has the strong Mankiewicz property, that is, every surjective isometry from $\mathcal{K}$ onto an arbitrary convex subset $L$ in a normed space $Y$ is affine.
\end{theorem}

We shall recall in the paragraphs below the notion of unitary element in unital JB$^*$-algebras. For the moment, we note that a Russo--Dye theorem, affirming that the closed convex hull of the unitary elements in a unital JB$^*$-algebra coincide with the closed unit ball, has been established by several authors. Holomorphic theory based on the tools introduced by Harris in \cite{Harris72} were employed by M.A. Youngson and J.D.M. Wright to counteract the lacking of polar decomposition (see \cite{WriYou78}). Another holomorphic tool, also due to Harris, is employed in the recent version of the Russo--Dye theorem for JB$^*$-triples proved by M. Mackey and P. Mellon in \cite[Theorem 3.1]{MackMell21}. A.A. Siddiqui obtained a new proof of the Wright--Youngson--Russo--Dye theorem for unital JB$^*$-algebra which does not rely on holomorphic theory, but on the fact that for each invertible element $a$ in a unital JB$^*$-algebra $\AAA$ there exists a unitary $u\in \AAA$ such that $a$ is positive and invertible in the $u$-isotope $\AAA(u)$ (cf. \cite{Sidd2010,Sidd2006,Sidd2007}, and see below for the concrete definition of the $u$-isotope). The reference \cite{Sidd2010} also contains a version of the results in \cite{KadPed} for unital JB$^*$-algebras. It follows from the just commented results that \begin{equation}\label{eq every convex body in a unital JB$^*$-algebra satisfies the strong Mankiewicz property} \begin{array}{l} \hbox{every convex body $K$ in a unital JB$^*$-algebra $\AAA$ satisfies the strong} \\
 \hbox{Mankiewicz property, that is, every surjective isometry from $K$ onto}  \\
 \hbox{an arbitrary convex subset of a normed space is affine.}
\end{array}
\end{equation}\smallskip

Our goal in this section is to prove a Jordan analogue of the commented result by Mori and Ozawa for maximal proper faces of the closed unit ball of a unital JB$^*$-algebra.\smallskip 

It is now an appropriate moment to recall the basic theory on unitaries in JB$^*$-algebras. We recall the notion of invertible element. An element $a$ in a unital JB$^*$-algebra $\mathfrak{A}$ is called \emph{invertible} if we can find $b\in \mathfrak{A}$ satisfying $a \circ b = \mathbf{1}$ and $a^2 \circ b = a.$ The element $b$ is unique and it will be denoted by $a^{-1}$ (cf. \cite[3.2.9]{HOS} or \cite[Definition 4.1.2]{CabGarPalVol1}). This definition agrees with the usual concept of invertibility in case that a unital C$^*$-algebra is regarded as a JB$^*$-algebra. The invertibility of an element $a$ in a JB$^*$-algebra $\AAA$ is equivalent to the invertibility of the $U_a$ operator in $B(\AAA)$, and in such a case $U_a^{-1} = U_{a^{-1}}$ (cf. \cite[Theorem 4.1.6]{CabGarPalVol1}). Actually, the notion of invertible elements works in the wider setting of unital Jordan Banach algebras, but in this paper we shall only work with unital JB$^*$-algebras. A procedure to produce invertible elements via exponentials reads as follows: the closed subalgebra of $\AAA$ generated by a fixed element $x\in \AAA$ and the unit is always associative, and hence the expression  $\displaystyle \exp(x):=\sum_{n=0}^{\infty}\frac{x^n}{n!}$ defines an element in $\AAA$ which is invertible with inverse $\exp(-x)$. The set of all invertible elements in $\AAA$ is an open set and the open unit ball of radius $1$ around the unit element is contained in the set of invertible elements \cite[Theorem 4.1.7]{CabGarPalVol1}.\smallskip
  
As in the case of unital C$^*$-algebra, an element $u$ in a unital JB$^*$-algebra $\AAA$ is called \emph{unitary} if it is invertible with inverse $u^*$. According to the standard notation, we shall write $\mathcal{U}(\AAA)$ for the set of all unitary elements in $\AAA$. This notion is consistent with the usual definition for C$^*$-algebras, that is, if a unital C$^*$-algebra is regarded as a JB$^*$-algebra both notions of unitaries coincide. From a wider point of view, unitary elements in a unital JB$^*$-algebra $\AAA$ fully coincide with the unitary tripotents in $\AAA$ when the latter is regarded as a JB$^*$-triple {\rm(}cf. \cite[Proposition 4.3]{BraKaUp78} or \cite[Theorem 4.2.24, Definition 4.2.25 and Fact 4.2.26]{CabGarPalVol1}{\rm)}. The reader should be warned that, even in the C$^*$-setting, the Jordan product of two unitaries in a unital JB$^*$-algebra might not be unitary. The set $\mathcal{U}(\AAA)$ is not, in general, connected. The connected component of $\mathcal{U}(\AAA)$ containing the unit element is called the \emph{principal component}, and it will be denoted by $\mathcal{U}^0(\AAA)$.\smallskip

We gather next some basic properties of the set of unitaries in a unital JB$^*$-algebra $\AAA$ which have been borrowed from several sources (they can be also found in \cite[Lemmata 2.1 and 2.2]{CuPe20TingleyUnitaries} and \cite{CueEnHirMiPer19}).
 \begin{enumerate}\item[$(\mathcal{PU}1)$]\label{PU1} For each $u\in \mathcal{U}(\AAA)$, the underlying Banach space of $\AAA$ becomes a unital JB$^*$-algebra with unit $u$ for the (Jordan) product defined by $x\circ_u y :=U_{x,y}(u^*)=\{x,u,y\}$ and the involution $*_u$ defined by $x^{*_u} :=U_u(x^*)=\{u,x,u\}$. (The JB$^*$-algebra $\AAA(u)= (\AAA,\circ_u,*_u)$ is called the \emph{$u$-isotope} of $\AAA$.) Furthermore $\mathcal{U}(\AAA) = \mathcal{U}(\AAA (u)) = \{\hbox{unitary tripotents in } \AAA\}$ \cite[Lemma 4.2.41]{CabGarPalVol1} and \cite[Proposition 4.3]{BraKaUp78}).
\item[$(\mathcal{PU}2)$]\label{PU2} There exist examples of unital JB$^*$-algebras $\AAA$ containing a unitary element $w$ for which no surjective linear isometry $T: \AAA \to \AAA$ maps $\mathbf{1}$ to $w$, that is, the group of all surjective linear isometries on $\AAA$ --equivalently, triple automorphisms-- on $M$ is not transitive on $\mathcal{U} (M)$ \cite[Example 5.7]{BraKaUp78}. It can be even shown that the unital JB$^*$-algebras $\AAA$ and $\AAA(w)$ are not Jordan $^*$-isomorphic \cite[Antitheorem 3.4.34]{CabGarPalVol1} and the connected component of $\mathcal{U} (\AAA)$ containing $w$ is not isometrically isomorphic to the principal connected component \cite[Remark 3.7]{CueEnHirMiPer19}.   
\item[$(\mathcal{PU}3)$]\label{PU3} The triple product of $\AAA$ satisfies $$\begin{aligned} \J xyz &= (x\circ y^*) \circ z + (z\circ y^*)\circ x -
		(x\circ z)\circ y^*\\
		&= (x\circ_{u} y^{*_{u}}) \circ_{u} z + (z\circ_{u} y^{*_{u}})\circ_{u} x -
		(x\circ_{u} z)\circ_{u} y^{*_{u}},
	\end{aligned}$$ for all $x,y,z\in \AAA$ and each $u\in \mathcal{U}(\AAA)$. That is, despite of the existence of many inequivalent unital JB$^*$-algebra structures on $\AAA$, there is only one triple product (follows from Kaup's Banach--Stone theorem \cite[Proposition 5.5]{Ka83}). 
\item[$(\mathcal{PU}4)$]\label{PU4} For each $u\in \mathcal{U}(\AAA)$, the mapping $U_u : \AAA\to \AAA$ is a surjective isometry and hence a triple isomorphism. Consequently, $U_u \left( \mathcal{U}(\AAA) \right) = \mathcal{U}(\AAA)$.  Furthermore, the operator $U_u : \AAA(u^*) \to \AAA(u)$ is a Jordan $^*$-isomorphism (cf. \cite[Theorem 4.2.28]{CabGarPalVol1}, \cite[Proposition 5.5]{Ka83} and \cite{WriYou78}).
\item[$(\mathcal{PU}5)$]\label{PU5} The principal connected component of $\mathcal{U} (\AAA)$ coincides with the intersection of the pricipal component of the set of invertible elements in $\AAA$ with $\mathcal{U} (\AAA)$ and admits the following geometric characterization:
\begin{equation}\label{eq algebraic charact principal component unitaries JBstar}\begin{aligned} \mathcal{U}^0(\AAA) &=\left\lbrace  U_{\exp(i h_n)}\cdots U_{\exp(i h_1)}(\textbf{1}) \colon n\in \mathbb{N}, \ h_j\in \AAA_{sa}, \ \forall\ 1\leq j \leq n  \right\rbrace \\
		&= \left\lbrace  u\in \mathcal{U} (\AAA) : \hbox{ there exists } w\in \mathcal{U}^0(\AAA) \hbox{ with } \|u-w\|<2 \right\rbrace
\end{aligned},\end{equation} and hence it is analytically arcwise connected \cite[Theorem 2.3]{CueEnHirMiPer19}.
\item[$(\mathcal{PU}6)$]\label{PU6} If $\AAA$ is a JBW$^*$-algebra, the set of all unitary elements in $\AAA$ is precisely the set given by the exponential of all skew-symmetric elements in $\AAA$, that is, $$\mathcal{U} (\AAA) = \big\{\exp(i h) : h\in \AAA_{sa}\big\} \ \ \  \hbox{(cf. \cite[Remark 3.2]{CuPe20TingleyUnitaries})}.$$ 
\end{enumerate}

Along this section we shall consider the natural projection  of $\mathbb{C}$ onto its unit sphere defined by $\mathrm{sgn}(\alpha)=\begin{cases}
	\frac{\alpha}{\|\alpha\|}  & \hbox{if $\alpha \ne 0$} \\
	0 & \hbox{ if $\alpha=0$}
\end{cases}$. We observe that for any JB$^*$-triple $E$, any element $x\in E$ and any $\alpha\in\CC$ the identity \begin{equation}\label{eq range alpha times x} \hbox{ $r(\alpha x)=\mathrm{sgn}(\alpha)r(x)$ holds}
\end{equation} where the range tripotent is computed in the second dual of $E$. Indeed, if we write $x^{\left[\frac{1}{3}\right]}=y$ ($n\in \mathbb{N}$), then $\left(\mathrm{sgn}(\alpha)|\alpha|^\frac{1}{3}y\right)^{[3]}=\mathrm{sgn}(\alpha)|\alpha|y^{[3]}=\alpha x$, so by the uniqueness of the cubic root $\left(\alpha x\right)^{\left[\frac{1}{3}\right]}=\mathrm{sgn}(\alpha)|\alpha|^\frac{1}{3}x^{\left[\frac{1}{3}\right]},$ and by induction $\left(\alpha x\right)^{\left[\frac{1}{3^n}\right]}=\mathrm{sgn}(\alpha)|\alpha|^\frac{1}{3^n}x^{\left[\frac{1}{3^n}\right]}$. For $\alpha=0$ everything is clear, for $\alpha\ne0$ we know that $|\alpha|^\frac{1}{3^n}\to 1$, so $r(\alpha x)=w^*-\lim_n\left(\alpha x\right)^{\left[\frac{1}{3^n}\right]}=w^*-\lim_n \mathrm{sgn}(\alpha)|\alpha|^\frac{1}{3^n}x^{\left[\frac{1}{3^n}\right]}= \mathrm{sgn}(\alpha)r(x)$.

\begin{lemma}\label{L_range_tripotent_is_principal_unitary} Let $\AAA$ be a unital JB$^*$-algebra and let $x\in S_\AAA$ and $\delta>0$. Then the range tripotent $r\left(\frac{(1+\delta) \mathbf{1}+x}{2+\delta}\right)$ lies in $\mathcal{U}^0(\AAA)$. 
\end{lemma}

\begin{proof} Denote $b\df \mathbf{1}+\frac{x}{1+\delta}$. Due to \eqref{eq range alpha times x} the equality  $r\left(\frac{(1+\delta)\mathbf{1}+x}{2+\delta}\right)=r\left(b\right)$ holds, in principle as elements in $\AAA^{**}$. Since $\|b-\mathbf{1}\|= \left\|\frac{x}{1+\delta}\right\|=\frac{1}{1+\delta}<1$, the element $b$ is invertible in $\AAA$. Therefore, by \cite[Remark 2.3]{JamPerSidTah2015} the range tripotent $r(b)$ is a unitary in $\AAA$. It is well known and  follows from the local Gelfand theory or from the triple functional calculus that $\|r(b)-b\|\le 1$. (Actually something stronger holds: since $b$ is invertible, it must be von Neumann regular in the sense in \cite{JamPerSidTah2015} and the triple spectrum of $b$ in $\AAA$ does not contain zero. Thus, the triple spectrum $\Omega_b$ is a subset of $[m,\|b\|]$ for some $m>0,$ and so $\|r(b)-b\|=\max\{ 1-m, \|b\|-1\} < 1$.) Furthermore, $\|b-\mathbf{1}\|=\left\|\frac{x}{1+\delta}\right\|=\frac{1}{1+\delta}<1$, so $\|r(b)-\mathbf{1}\|<2$. Therefore we deduce from \eqref{eq algebraic charact principal component unitaries JBstar}  (equivalently, \cite[Theorem 2.3]{CueEnHirMiPer19} or \cite[Lemma 2.2]{CuPe20TingleyUnitaries}) that $r(b)\in\mathcal{U}^0(\AAA)$.
\end{proof}

The next technical lemma plays a key role in our arguments. We shall avoid the difficulties produced by the lacking of polar decompositions for elements in JBW$^*$-algebras by employing isotopes. 

\begin{lemma}\label{L_construction_of_unitaries_in_F_p}
Let $\AAA$ be a unital JB$^*$-algebra and $F_e\subset S_\AAA$ be a norm-closed proper face of the closed unit ball of $\AAA$ {\rm(}associated with a compact tripotent $e\in\AAA^\dd${\rm)}. Then for any $w\in{\mathcal U}^0(\AAA)\cap F_e$, any $x\in F_e$, and any $\eps>0,$ there are unitaries $u,v\in F_e\cap {\mathcal U}^0(\AAA)$ such that $$\|u+v-(x+w)\|<\eps.$$
\end{lemma}

\begin{proof} Fix an arbitrary  $0< \delta<1$ and denote $a\df\frac{(1+\delta)w+x}{2+\delta}$ and $r\df r(a)$. Lemma~\ref{L_range_tripotent_is_principal_unitary} implies that $r$ lies in the principal connected component of the $w$-isotope $\AAA (w)$, which is denoted by $\mathcal{U}^0(\AAA(w))$. This means, that $r$ lies in the same connected component of $\mathcal{U}(\AAA)$ as $w$ does. Therefore, since $w\in\mathcal{U}^0(\AAA)$, it follows that $r\in{\mathcal U}^0(\AAA)$. \smallskip
	
We know that $a\ge0$ in the $r$-isotope $\AAA(r),$ and since $$\|a^2\|=\left\| \left(\frac{(1+\delta)w+x}{2+\delta}\right)^2\right\|\le \left\| \frac{(1+\delta)w+x}{2+\delta}\right\|^2\le \left(\frac{(1+\delta)+1}{2+\delta}\right)^2=1,$$ where all squares are computed in the $r$-isotope $\AAA(r)$, we also have $0\leq a^2 = \{a,r,a\} \le r$ (i.e. we consider the square $a^2$ and the order in the isotope $\AAA(r)$). We define $u,v\df a\pm i\sqrt{r-a^2}$ (again, all operations are considered in the isotope $\AAA(r)$). Then $u,v$ are unitaries in the principal connected component of the isotope $\mathcal{U}(\AAA(r))(=\mathcal{U}(\AAA)),$ that is the connected component of $\mathcal{U}(\AAA))$ containing $r,$ because $\|r-u\|,\|r-v\|\leq \sqrt{2}$ (see \eqref{eq algebraic charact principal component unitaries JBstar} or \cite[Theorem 2.3]{CueEnHirMiPer19}). So, both $u,v$ lie in the principal connected component ${\mathcal U}^0(\AAA)$.\smallskip

Now, we deduce from the identity $\frac{u+v}{2}=\frac{1+\delta}{2+\delta}w+\frac{1}{2+\delta}x =a\in F_e,$ combined with the fact that $F_e$ is a  a proper closed face of the closed unit ball, that $u,v\in F_e$.\smallskip
 
Finally, we check the desired distance $$\|v+u-x-w\| = \|2 a -x-w\|=\left\|2\frac{(1+\delta)w+x}{2+\delta}-w-x\right\|=\left\|\frac{\delta w-\delta x}{2+\delta}\right\|\le \frac{2\delta}{2+\delta},$$ witnessing that for $\delta$  small enough (specifically for $\delta<\frac{2\eps}{2-\eps}$) we get $\|v+u-x-w\|<\eps$, as desired.
\end{proof}

Let us observe that the compact tripotent $e\in \AAA^{**}$ in the statement of Lemma~\ref{L_construction_of_unitaries_in_F_p} played no real role in the arguments, which are more geometric. \smallskip

We shall conclude this section with the pursued Jordan version of the result proved by Mori and Ozawa for unital C$^*$-algebras in \cite[Lemma 18]{MoriOza2018}.  

\begin{theorem}\label{Russo_Dye_theorem_for_F_p}
Let $\AAA$ be a unital JB$^*$-algebra and let $F_e\subset S_\AAA$ be a norm-closed proper face associated with a compact tripotent $e\in\AAA^\dd$. Suppose $u$ is a unitary in $\AAA$ satisfying that $e$ is a projection in the JBW$^*$-algebra $\AAA^{**}(u)$ given by the $u$-isotope of $\AAA^{**}$. Let ${\mathcal U}^0(\AAA (u))$ denote the principal connected component of the set of unitary elements in $\AAA (u)$. Then $F_e=\overline{\mathrm{co}}\left(F_e\cap{\mathcal U}^0(\AAA(u))\right)$. The conclusion clearly holds for every norm-closed proper face of $\mathcal{B}_{\AAA}$ associated with a compact projection in $\AAA^{**}$ and ${\mathcal U}^0(\AAA)$. 
\end{theorem}

\begin{proof} Let us first observe that, since $u\in \AAA$, by weak$^*$ density of $\AAA$ in its bidual, the bidual of $\AAA(u),$ $\left(\AAA(u)\right)^{**}$, identifies with the JBW$^*$-algebra  $\AAA^{**}(u)$. It is not hard to check that $F$ is a closed proper face of the closed unit ball of $\AAA(u)$ associated with a compact projection in $\left(\AAA(u)\right)^{**}$. So, up to replacing the original Jordan product and involution by $\circ_u$ and $^{*_u}$, we can always assume that $u$ is the unit element in $\AAA,$ $e=p$ is a projection in $\AAA^{**}$ and ${\mathcal U}^0(\AAA (u)) ={\mathcal U}^0(\AAA)$ is the principal connected component of $\mathcal{U} (\AAA)$.\smallskip
	  
The inclusion $\supseteq$ is obvious.\smallskip 
	
Take any $x\in F_p$ and fix a positive $0<\eps <1$. We shall prove the existence of two sequences $(u_n)_n,$ and $(v_n)_n$ in $\mathcal{U}^0 (\AAA)\cap F_p$ such that \begin{equation}\label{eq first statement by induction} \|u_n+v_n-x-v_{n-1}\|< \frac{\eps}{3}.
\end{equation} Set $v_0\df\mathbf{1}$, which clearly lies in both sets ${\mathcal U}^0(\AAA)$ and $F_p$. By applying Lemma~\ref{L_construction_of_unitaries_in_F_p} with $w =v_0 = \mathbf{1}\in{\mathcal U}^0(\AAA)\cap F_p$, $x\in F_p$, and $\frac{\eps}{3}>0,$ there are unitaries $u_1,v_1\in F_p\cap {\mathcal U}^0(\AAA)$ such that $$\|u_1+v_1-(x+v_0)\|<\frac{\eps}3.$$ Suppose, we have already defined $u_1,\ldots, u_n$ and $v_1,\ldots, v_n$ satisfying \eqref{eq first statement by induction}, a new application of Lemma~\ref{L_construction_of_unitaries_in_F_p} with $w =v_n \in{\mathcal U}^0(\AAA)\cap F_p$, $x\in F_p$, and $\frac{\eps}{3}>0,$ gives the elements $u_{n+1}$ and $v_{n+1}$ in $F_p\cap {\mathcal U}^0(\AAA)$ satisfying the desired property. \smallskip

Now, by an easy induction argument we deduce from \eqref{eq first statement by induction} that  $$\left\|u_1+u_2+\dots+u_n+v_n-\mathbf{1}-n x\right\|< n \frac{\eps}{3}$$ for all natural $n$. Therefore $$\left\|\frac{u_1+u_2+\dots+u_n}{n}+\frac{v_n-\mathbf{1}}{n} - x\right\| < \frac{\eps}{3} \ \ \ \ \hbox{($\forall n\in \mathbb{N}$)}.$$ Take any natural $k$ such that $k\ge\frac{3}{\eps}$. Then $$\begin{aligned} \left\|\frac{u_1+u_2+\dots+u_k}{k}-x\right\|&\le\left\|\frac{u_1+u_2+\dots+u_k}{k}+\frac{v_k-\mathbf{1}}{k}-x\right\|+\frac{\|v_k-\mathbf{1}\|}{k} \\ &<\frac{\eps}{3}+2\frac{\eps}{3}=\eps.
	\end{aligned}$$ Since $\frac{u_1+u_2+\dots+u_k}{k}\in\mathrm{co}\left(F_p\cap{\mathcal U}^0(\AAA)\right)$, we conclude that $x\in \overline{\mathrm{co}}\left(F_p\cap{\mathcal U}^0(\AAA)\right)$.
\end{proof}

An strengthened version of Theorem~\ref{Russo_Dye_theorem_for_F_p} valid for all maximal proper faces and the principal connected component of $\mathcal{U} (\AAA)$ will be obtained at the end of section~\ref{sec: domination of minimal tripotents in the second dual by unitaries in the main component} (see Theorem \ref{Russo_Dye_theorem_for_maximal proper faces}). At this moment it is not clear why each minimal tripotent in the second dual of a unital JB$^*$-algebra $\AAA$ is dominated by a unitary in the principal component of $\mathcal{U} (\AAA)$. \smallskip

In this note we shall also need to widen the known list of  Banach spaces satisfying the strong Mankiewicz property. In \cite[Lemma 19]{MoriOza2018} it if shown that if $\psi: A\to \mathbb{C}$ is a  non-zero multiplicative linear functional on a unital C$^*$-algebra, then the closed unit ball of the real C$^*$-algebra $A 
_{\mathbb{R}}^{\psi}:=  \psi^{-1}(\mathbb{R})$  coincides with the closed convex hull of the unitary elements in $A 
_{\mathbb{R}}^{\psi}$, and hence it satisfies the strong Mankiewicz property. The arguments in \cite{MoriOza2018} rely on a Russo--Dye theorem for real von Neumann algebras (\cite[Theorem 7.2.4]{Li2003}, see also \cite[proof of Corollary 3]{MoriOza2018}). In the next lemma we shall pursue a version of the just commented result for JB$^*$-algebras. The main handicap is the lack of an appropriate Russo--Dye theorem for unital real JBW$^*$-algebras. That remains as a deep open problem. For our purposes we shall employ an alternative trick.  \smallskip

We recall first that a \emph{real JB$^*$-algebra} is a closed $^*$-invariant real subalgebra of a (complex) JB$^*$-algebra. Unital real JB$^*$-algebras were introduced by K. Alvermann in \cite{Alver86}. Besides the just quoted reference an axiomatic intrinsic definition of unital real JB$^*$-algebras can be found in \cite{PeAxiom2003}. The class of real JB$^*$-algebras includes all JB-algebras, which are precisely those real JB$^*$-algebras whose involution is the identity mapping, and all real C$^*$-algebras (i.e., closed $^*$-invariant real subalgebras of C$^*$-algebras) equipped with their natural Jordan product \cite{Li2003,GoodearlBook82}. \smallskip 

Real JB$^*$-algebras can be also understood as examples of real JB$^*$-triples \cite[p. 321]{IsKaRo95}. A \emph{real JB$^*$-triple} is by definition a real closed subtriple of a JB$^*$-triple (see \cite{IsKaRo95}). Every JB$^*$-triple is a real JB$^*$-triple when it is regarded as a real Banach space. As in the case of real C$^*$-algebras, real JB$^*$-triples can be obtained as \emph{real forms} of JB$^*$-triples. More concretely, given a real JB$^*$-triple $E$, there exists a unique (complex) JB$^*$-triple structure on its algebraic complexification $X= E \oplus i E,$ and a conjugation (i.e. a conjugate linear isometry of period 2) $\tau$ on $X$ such that $$E = X^{\tau} = \{ z\in X : \tau (z) = z\},$$ (see \cite{IsKaRo95}). Consequently, every real C$^*$-algebra is a real JB$^*$-triple with respect to the product given in \eqref{eq triple product JCstar triple}, and the Banach space $B(H_1,H_2)$ of all bounded real linear operators between two real, complex, or quaternionic Hilbert spaces also is a real JB$^*$-triple with the same triple product.\smallskip 

The real C$^*$-algebra $C([0,1],\mathbb{R})$ of all real-valued continuous functions on $[0,1]$ also is a real JB$^*$-algebra for which the convex hull of its unitary elements reduces to the constant functions. So, it is hopeless to expect a Russo--Dye theorem for real JB$^*$-algebras, the closest statement can be found in \cite[Theorem 7.2.4]{Li2003}. As we commented before, the conclusion holds under additional hypotheses. A real von Neumann algebra is a real C$^*$-algebra which is also a dual Banach space. Similarly, a real JBW$^*$-triple (respectively, a real JBW$^*$-algebra) is a real JB$^*$-algebra which is a dual Banach space. Each real JBW$^*$-triple admits a unique isometric predual and its triple product is separately weak$^*$ continuous \cite{MarPe00}. The bidual of each real JB$^*$-triple is a real JBW$^*$-triple \cite[Lemma 4.2]{IsKaRo95}.\smallskip

\begin{proposition}\label{p strong Mankiewicz for the real JB$^*$-algebras associated with a Jordan multiplicative functional} Let $\AAA$ be a unital JB$^*$-algebra, and suppose that $\psi: \AAA\to \mathbb{C}$ is a non-zero Jordan $^*$-homomorphism. Consider the unital real JB$^*$-algebra $\AAA_{\mathbb{R}}^{\psi}:=  \psi^{-1}(\mathbb{R})$. Then $\frac12 \mathcal{B}_{\AAA_{\mathbb{R}}^{\psi}}$ is contained in the closed convex hull of the unitary elements in $\AAA_{\mathbb{R}}^{\psi}$, and consequently every convex body $K \subseteq \AAA_{\mathbb{R}}^{\psi}$ satisfies the strong Mankiewicz property.
\end{proposition}

\begin{proof} Let us recall a device from \cite[Lemma 7.2.5]{Li2003} which remains valid for all unital real JB$^*$-algebras. Let $x$ be an element in the closed unit ball of a unital real JB$^*$-algebra $\mathcal{A}$. Let us write $x = h + k$ with $h = h^*$ and $k = -k^*$. The unital real JB$^*$-subalgebras $\mathcal{A}_{h,\mathbf{1}}$ and $\mathcal{A}_{k,\mathbf{1}}$ of $\mathcal{A}$ generated by $\{h,\mathbf{1}\}$ and $\{k,\mathbf{1}\}$, respectively, are isometrically JB$^*$-algebra isomorphic to unital abelian real C$^*$-algebras, and hence we can apply \cite[Lemma 7.2.5]{Li2003} and its proof to find $b=b^*\in  \mathcal{A}_{h,\mathbf{1}}$ and $a\in \mathcal{A}_{k,\mathbf{1}}$ such that $ h = \cos(b)$ and $\displaystyle k = \frac{\exp(a) - \exp(-a)}{2}$. Therefore $\displaystyle \frac{x}{2} = \frac12 \cos(b) + \frac14 \exp(a) + \frac14 (- \exp(-a))$. That is, every element in $\frac12 \mathcal{B}_{\mathcal{A}}$ is a convex combination of three elements of the form $\cos(b)$, $\exp(a)$ and $\exp(-a)$ with $b = b^*$ and $a = -a^*$ in $\mathcal{A}$.\smallskip
	
If we additionally assume that $\mathcal{A}$ is a real JBW$^*$-algebra, and the real JBW$^*$-subalgebra generated by $b$ and $\mathbf{1}$ is isometrically isomorphic as real JBW-algebra to some $C(K, \mathbb{R})$ for an appropriate hyperStonean space $K$ (cf. \cite[Lemma 4.1.11]{HOS} or \cite[\S 6.3]{Li2003}). It is well known that the convex hull of the unitaries in $C(K,\RR)$ is dense in $\mathcal{B}_{C(K,\RR)}$ if and only if $K$ is totally disconnected (see \cite{Bade57}), therefore $\cos(b)$ is contained in the closed convex hull of the unitaries in the real von Neumann algebra generated by $b$, which are clearly unitaries in $\mathcal{A}$. We have shown the following: \begin{equation}\label{eq the ball of radius one half is in the convex hull of unitaries real}\begin{array}{l} \hbox{For each real JBW$^*$-algebra $\mathcal{A}$ the set $\frac12 \mathcal{B}_{\mathcal{A}}$ is contained}\\ 
		 \hbox{in the closed convex hull of the unitary elements in $\mathcal{A}$.}
	\end{array}  
\end{equation}

Let us return to the original setting of this proposition. We can clearly regard $\psi $ as weak$^*$ continuous Jordan $^*$-homomorphism from $\AAA^{**}$ to $\mathbb{C}$.  Clearly, the space $\left(\AAA^{**}\right)_{\mathbb{R}}^{\psi}:= \{a\in \AAA^{**} : \psi (a)\in \mathbb{R}\}$ is a real JBW$^*$-algebra of $\AAA^{**}$. Let us pick an element $a$ in the real JB$^*$-algebra $\AAA_{\mathbb{R}}^{\psi}\subseteq \left(\AAA^{**}\right)_{\mathbb{R}}^{\psi}$ with $\|a\|\leq \frac12$, a non-zero $\phi\in {\AAA^{*}}$ and $\varepsilon>0$. By \eqref{eq the ball of radius one half is in the convex hull of unitaries real}, there exist unitary elements $u_1,\ldots,u_m\in \left(\AAA^\dd\right)_{\mathbb{R}}^{\psi}$ and $t_1,\ldots,t_m\in (0,1)$ satisfying $\sum_j t_j =1$ and  $\left\|a-\sum_j t_j u_j\right\|<\frac{\varepsilon}{3 \|\phi\|}.$ When the unitaries $u_1,\ldots,u_m$ are regarded in $\AAA^{**}$ we can find self-adjoint elements $h_1,\ldots, h_m$ in $\AAA^{**}$ such that $u_j = \exp(i h_j)$ for all $j=1,\ldots, m$. Since $u_j \in \left(\AAA^\dd\right)_{\mathbb{R}}^{\psi}$ we deduce that $\psi (u_j)\in \mathbb{R},$ and thus $\psi (h_j)\in \pi \mathbb{Z},$ for all $j=1,\ldots, m$. \smallskip

The strong$^*$-topology of $\AAA^{**}$ in the sense of \cite[\S 3]{BarFri90}, when regarded as a JBW$^*$-triple, is a topology stronger than the weak$^*$-topology, makes the Jordan and the triple product of $\AAA^{**}$ jointly strong$^*$ continuous on bounded sets \cite{RodPa91,PeRo2001}, and an appropriate \emph{Kaplansky density theorem} holds for this topology \cite[Corollary 3.3]{BarFri90}. That is, $\mathcal{B}_{\AAA}$ is strong$^*$-dense in $\mathcal{B}_{\AAA^{**}}$. Furthermore, by the mentioned properties, $\mathcal{B}_{\AAA_{sa}}$ is strong$^*$-dense in $\mathcal{B}_{\AAA^{**}_{sa}}$. So, for each $j=1,\ldots,m$, we can find a bounded net $(h^{j}_{\lambda_j})_{\lambda_j}$ in $\AAA_{sa}$ converging to $h_j$ with respect to the strong$^*$-topology. So, we can find $h_0^j\in \AAA_{sa}$ for which $$\left|\phi\left( \exp(i h_{0}^{j})  - \exp(i h_j)  \right)\right|<\frac{\varepsilon}{3}, \hbox{ and } \left| 1-\exp\left(-i \psi\left( h_{0}^{j}  -  h_j \right)\right)\right|<\frac{\varepsilon}{3\|\phi\|},$$ for all $j = 1,\ldots, m.$ Let us define $\hat{h}_{0}^{j} = {h}_{0}^{j} - \left(\psi\left( {h}_{0}^{j}\right) - \psi\left(h_j \right) \right)\mathbf{1}\in \AAA_{sa}$. We know by definition that $\psi \left(\hat{h}_{0}^{j} \right) = \psi \left({h}_{j} \right)\in \pi\mathbb{Z},$ and by the properties of $\psi$, $\psi\left( \exp(i\hat{h}_{0}^{j}) \right) =\exp\left( i\psi\left( \hat{h}_{0}^{j} \right) \right)\in \mathbb{R}.$  Therefore, the element $\exp(i\hat{h}_{0}^{j})$ is a unitary in $\AAA_{\mathbb{R}}^{\psi}$ for all $j=1,\ldots, m$. It follows by construction that $$\begin{aligned} \left|\phi \left( a-\sum_j t_j \exp(i \hat{h}_{0}^{j}) \right) \right| &\leq \left|\phi \left( a-\sum_j t_j \exp(i {h}_{j}) \right) \right| \\ &+ \left|\phi \left(  \sum_j t_j \left( \exp(i {h}_{j}) - \exp(i {h}_{0}^{j})\right)  \right) \right| \\
& +  \left|\phi \left( \sum_j t_j \left( \exp(i {h}_{0}^{j}) - \exp(i \hat{h}_{0}^{j})\right)  \right) \right| \\
&\leq  \|\phi\| \left\|a-\sum_j t_j u_j\right\| +  \sum_j t_j \left|\phi \left( \exp(i {h}_{j}) - \exp(i {h}_{0}^{j})\right)\right| \\
&+ \|\phi\|  \sum_j t_j \left\| \exp(i {h}^j_{0}) - \exp( i \hat{h}_{0}^{j}) \right\| \\
< 2 \|\phi\| \frac{\varepsilon}{3\|\phi\|}  + \|\phi\| &  \sum_j  t_j \left\| \exp(i {h}^j_{0})\circ \left(\mathbf{1}- \exp\left(-i \left(\psi\left( {h}_{0}^{j}\right) - \psi\left(h_j \right) \right)\mathbf{1}\right)\right)\right\| \\
\leq 2 \frac{\varepsilon}{3}  + \|\phi\|  \sum_j& t_j  \left\| \exp(i {h}^j_{0}) \right\| \left| {1}- \exp\left(-i \left(\psi\left( {h}_{0}^{j}\right) - \psi\left(h_j \right) \right)\right)\right| <\varepsilon.
\end{aligned} $$ We have therefore shown that every $a\in \frac12 \mathcal{B}_{\AAA_{\mathbb{R}}^{\psi}}$ lies in the weak-closure of the convex hull of $\mathcal{U} \left(\AAA_{\mathbb{R}}^{\psi}\right)$, which by the Hahn--Banach theorem is precisely the norm convex hull of $\mathcal{U} \left(\AAA_{\mathbb{R}}^{\psi}\right)$. Theorem 2 in \cite{MoriOza2018} gives the final conclusion. 	
\end{proof}

\section{Affine behaviour of the isometry on the maximal proper faces of the closed unit ball}\label{sec: affine behaviour on max faces}

This section is devoted to prove the third main tool required in our arguments. The main goal will consist in proving that every surjective isometry from the unit sphere of a unital JB$^*$-algebra $\AAA$ onto the unit sphere of any other Banach space is affine on every maximal proper face determined by a minimal projection in $\AAA^{***}$, a result proved by Mori and Ozawa in the case of unital C$^*$-algebras \cite[Proposition 20]{MoriOza2018}. To make easier the arguments we shall establish first a series of technical lemmata. In this section our conclusion will be limited to maximal proper faces associated with a minimal projection in the second dual, the conclusion for general minimal tripotent in the second dual will require another geometric tool developed in the next section.

\begin{lemma}\label{L_Existence_of_suitable_approximate_units} Let $L$ be a {\rm(}non-necessarily unital{\rm)} JB$^*$-subalgebra of a unital JB$^*$-algebra $\AAA$. Then there exist a directed set $\Lambda$ and nets $(d_\lambda)_{\lambda\in\Lambda},(e_\lambda)_{\lambda\in\Lambda}, (f_\lambda)_{\lambda\in\Lambda}\subset L$ such that for any $\lambda\in\Lambda$ the elements $d_\lambda,$ $e_\lambda$ and $f_\lambda$ belong to a commutative and associative JB$^*$-subalgebra of $L$ and the following statements hold: $$0\le e_\lambda\le d_\lambda\le f_\lambda\le \mathbf{1},\ \ e_\lambda\circ f_\lambda=e_\lambda, \ \  d_\lambda =d_\lambda\circ f_\lambda, \hbox{ and } e_\lambda=d_\lambda^2.$$ Furthermore,  for each $x\in L$ the nets $(d_\lambda\circ x)_{\lambda\in\Lambda}$, $(e_\lambda\circ x)_{\lambda\in\Lambda}$ and $(U_{d_\lambda}(x))_{\lambda\in\Lambda}$ all converge in norm to $x$.
\end{lemma}

\begin{proof} It is known that each JB$^*$-algebra admits an approximate unit (cf. \cite[Lemma 4 and the preceding comments]{Ed80} or \cite[Proposition 3.5.23]{CabGarPalVol1}), so we can find a directed set $\Gamma$ and a net $(i_\kappa)_{\kappa\in \Gamma}\subset L$ of positive elements bounded by $\mathbf{1}$ which is an approximate unit in $L$ and satisfying that $(i_\kappa^2)_{\kappa\in \Gamma}$ also is an approximate unit in $L$. Now for a fixed $i_\kappa$ we consider the JC$^*$-subalgebra $B_{i_\kappa}$ of $L$ generated by $i_\kappa$. Since $i_\kappa$ is self-adjoint, $\B_{i_\kappa}$ is associative and therefore it is isometrically isomorphic to a commutative C$^*$-algebra. Since $\|i_\kappa\| =1$ and $0\le i_\kappa\le \mathbf{1}$, $\B_{i_\kappa}$ is isometrically isomorphic to $C_0(K)$ for some compact Hausdorff space $K \subseteq [0,1]$ with $1\in K,$ and under this identification $i_\kappa$ corresponds to the embedding of $K$ into $\mathbb{C}$. \smallskip
	
For any $n\in\N$, $n>1$ let $d_{\kappa,n}$  (respectively, $f_{\kappa,n}$) be the element of $B_{i_\kappa}$ associated with the continuous function in $C_0(K)$ defined by $$t\mapsto\chi_{\left[\frac{1}{n},1\right]}(t) \frac{n t-1}{n-1},$$ $\big($respectively, $t\mapsto n\chi_{\left[0,\frac{1}{n}\right]}(t) \ t+\chi_{\left(\frac{1}{n},1\right]}(t)\big)$. Finally, let $e_{\kappa,n}\df d_{\kappa,n}^2$ and set $\Lambda\df \Gamma\times\N.$ Since $d_{\kappa,n}, f_{\kappa,n}$ are both associated with functions vanishing at $0$, they lie in $\B_{i_\kappa}\subset L$, and therefore $e_{\kappa,n}\in L$ as well. \smallskip
	
By construction we already have $e_{\lambda}=d_{\lambda}^2$. Furthermore $0\le e_\lambda \le d_\lambda\le f_\lambda\le \mathbf{1}$ follows from the fact that the same identities obviously hold for the associated functions in $C([0,1])$. By the same reasons the identities $e_\lambda\circ f_\lambda=e_\lambda$ and $d_\lambda=d_\lambda\circ f_\lambda$  are also true.\smallskip
	
Since for fixed $\kappa\in \Gamma$ the sequence $d_{\kappa,n}$ (respectively, $e_{\kappa,n}$) is associated with an increasing sequence of continuous functions converging pointwise to the function associated with $i_\kappa$ (respectively, with $i_\kappa^2$) which also continuous, this convergence is uniform by Dini's theorem, so $d_{\kappa,n}$ (respectively, $e_{\kappa,n}$) converges to $i_\kappa$ (respectively, $i_\kappa^2$) in norm. Notice that these convergences are uniform in $\kappa$ due to the uniform convergence of the associated functions on the entire interval $[0,1]$.\smallskip
	
Let us fix a non-zero $x\in L$ and $\eps>0$. For any $\lambda=(\kappa_0,n_0)\in \Gamma\times \N$ we can find $\kappa_1>\kappa_0$ such that $\|i_\kappa\circ x-x\|<\frac{\eps}{2}$ and $\|i_\kappa^2\circ x-x\|<\frac{\eps}{2}$ for all $\kappa\ge\kappa_1$. We can also find $n_1\ge n_0$ such that $\|d_{\kappa,n}-i_{\kappa}\|<\frac{\eps}{2\|x\|}$ and $\|e_{\kappa,n}-i_{\kappa}^2\|<\frac{\eps}{2\|x\|}$ for all $\kappa\in\Gamma$ and $n\ge n_1$. Then $$\|d_{\kappa,n}\circ x-x\|\le\|d_{\kappa,n}\circ x-i_\kappa\circ x\|+\|i_{\kappa}\circ x-x\|< \|d_{\kappa,n}-i_\kappa\|\|x\|+\frac{\eps}{2}<\eps,$$ and $$\|e_{\kappa,n}\circ x-x\|\le\|e_{\kappa,n}\circ x-i^2_\kappa\circ x\|+\|i_{\kappa}^2\circ x-x\|< \|e_{\kappa,n}-i_\kappa^2\|\|x\|+\frac{\eps}{2}<\eps$$ for all $n\ge n_1$, $\kappa\ge\kappa_1$.\smallskip 

It only remains to prove that $U_{d_\lambda}(x)\to x$ in norm. Since $U_{d_\lambda}(x)=2d_\lambda\circ(d_\lambda\circ x)-e_\lambda \circ x$, all we need to prove is that $d_\lambda\circ(d_\lambda\circ x)\to x$ in norm. For this purpose fix $\lambda_0\in\Lambda$ and $\eps>0$. Then we take $\lambda>\lambda_0$ such that $\|d_\lambda\circ x-x\|\le\frac{\eps}{2}$. Then from $\|d_\lambda\|\le1$ we get $$\|d_\lambda\circ(d_\lambda\circ x)-x\|\le\|d_\lambda\circ(d_\lambda\circ x-x)\|+\|d_\lambda\circ x-x\|< \frac{\eps}{2}\|d_\lambda\|+\frac{\eps}{2}\le \eps.$$  
\end{proof}

\begin{lemma}\label{L_Jordan_product_estimate_for_selfadjoint_elements} Let $A$ be a C$^*$-algebra, let $p,q$ be projections in $A$ with $0\le q\le p$. Suppose $h$ is an element in $A_{sa}$, and $\eps>0$. If $\|p\circ h\|\leq \eps$, then $\|q\circ h\|\leq 3\eps$.
\end{lemma}

\begin{proof} Let $h_i\df P_i(p)(h)$. Then $h_2=php$, $h_1=ph(\jed-p)+(\jed-p)hp$, and $p\circ h=h_2+\frac{1}{2}h_1$. Since Peirce projections are contractive, we get $\|h_2\|\leq\eps$, $\|h_1\| \leq 2\eps$. Since $hp=php+(\jed-p)hp=h_2+h_1p$, we obtain $\|hp\|\le\|h_2\|+\|h_1\|\|p\|\leq 3\eps$.\smallskip
	
Having in mind that the mapping $U_h$ is positive ($h$ is self-adjoint), the assumption $0\le q\le p$, gives $0\le hqh\le hph$. Therefore, we deduce from the Gelfand--Naimark axiom that  $$\|qh\|^2=\|hq\|^2=\|hqh\|\le\|hph\|=\|hp\|^2\leq (3\eps)^2,$$ which leads to $\|q\circ h\|\le\frac{\|qh\|+\|hq\|}{2}\leq 3\eps$.
\end{proof}

\begin{lemma}\label{L_Jordan_product_estimate_for_general_elements} Let $p,q$ be projections in a JB$^*$-algebra $\AAA$ with $0\le q\le p$. Suppose $x$ is an element in $\AAA$ satisfying $\|p\circ x\|\leq\eps$ for some $\eps>0$. Then $\|q\circ x\|\leq 6 \eps$.
\end{lemma}

\begin{proof} Let $x = h+ i k$ with $h,k\in \AAA_{sa}$.  We note that $p\circ x=p\circ h+ip\circ k$ corresponds to the decomposition of $p\circ x$ in $\AAA_{sa}\oplus i \AAA_{sa}$. Therefore $\|p\circ h\|,\|p\circ k\|\le\|p\circ x\|\leq \eps$.\smallskip
	
The element $a=(p-q) -q$ is self-adjoint and hence the JB$^*$-subalgebra $\mathfrak{B}$ of $\AAA$ generated by $h$ and $a$ is a JC$^*$-algebra by the Shirshov-Cohn theorem \cite[Theorem 2.4.14]{HOS}. The element $a^2 =  (p-q) +q$ lies in $\mathfrak{B},$ and thus $p-q$, $q$ and $h$ lie in $\mathfrak{B}$.  Lemma~\ref{L_Jordan_product_estimate_for_selfadjoint_elements} proves that $\|q\circ h\|\leq 3\eps$. We similarly obtain that $\|q\circ k\|\leq 3\eps$ and hence $\|q\circ x\|\le \|q\circ h\|+\|q\circ k\|\le 6\eps.$
\end{proof}

Our main result in this section will show that every surjective isometry from the unit sphere of a unital JB$^*$-algebra $\AAA$ onto the unit sphere of any Banach space is affine when restricted to any maximal proper face of $\mathcal{B}_{\AAA}$ associated with a minimal projection in $\AAA^{**}$. The arguments here are inspired by a similar result proved by Mori and Ozawa in the setting of C$^*$-algebras (cf. \cite[Proposition 20]{MoriOza2018}). There is no need to say that the setting of JB$^*$-algebras requires a non-trivial adaptation with an appropriate combination of Jordan techniques.

\begin{proposition}\label{Prop_20_from_Mori-Ozawa}
	Let $\AAA$ be a unital JB$^*$-algebra, $Y$ a Banach space, $\Delta:S_\AAA\to S_Y$ a surjective isometry, and $F_p\subset S_{\AAA}$ a maximal norm-closed proper face associated with a minimal projection $p\in \AAA^\dd$. Then there exist a directed set $\Lambda$ and a net $(\Theta_\lambda)_\lambda$ of affine contractions $\Theta_\lambda: F_p\to K_\lambda,$ where each $K_\lambda$ is a closed convex subset of $ F_p$, such that $\Theta_\lambda\to\mathrm{id}_{F_p}$ in the point-norm topology {\rm(}i.e. in the strong operator topology{\rm)}, each $\Delta(K_\lambda)$ is convex, and each $K_\lambda$ is affinely isometrically isomorphic to the closed ball of a real JB$^*$-algebra $\mathcal{R}_\lambda$ such that $\B_{\mathcal{R}_\lambda}$ satisfies the strong Mankiewicz property {\rm(}and therefore $\Delta\zuz_{F_p}$ is affine{\rm )}.  
\end{proposition}

\begin{proof} Let $\phi\in \AAA^\ast$ be the pure state associated with $p$. Let us consider the standard seminorm $\|a\|_\phi^2\df\phi(a\circ a^\ast),$ and set $L\df \{a\in \AAA:\|a\|_\phi=0\}$.\smallskip
	
Since $\|a\|_\phi=\|a^\ast\|_\phi$ for all $a\in \AAA$, $L$ is clearly self-adjoint. Furthermore,  $p$ is minimal implies that it is compact \cite[Theorem 3.4]{BuFerMarPe}, so $\mathbf{1}-p$ is open and so $$L=\AAA_0^\dd(p)\cap \AAA=\AAA_2^\dd(\mathbf{1}-p)\cap \AAA$$ (see, for example, any of the references \cite[Lemma 1.5]{FriRu85}, \cite{Pe15}, \cite[Proposition 1.2]{BarFri87Groth}, \cite[Lemma 4.1]{EdRu88} or \cite[Lemma 6.11$(e)$]{HKPP-MWnK} for the first equality). The identity $L=\AAA_2^\dd(\mathbf{1}-p)\cap \AAA$ implies that $L$ is a closed inner inner ideal of $\AAA$. Consequently, for any $a,b\in L$ we have $$a\circ b=\{a,\mathbf{1},b\}=\{a,\mathbf{1}-p,b\}\in \AAA_0^\dd(p).$$ Therefore $L$ is a Jordan $^*$-subalgebra of $\AAA$.\smallskip 
	
By Lemma~\ref{L_Existence_of_suitable_approximate_units} we can find a directed set $\Lambda$ and nets $(d_\lambda)_{\lambda\in\Lambda},(e_\lambda)_{\lambda\in\Lambda}, (f_\lambda)_{\lambda\in\Lambda}\subset L= \AAA_0^\dd(p)\cap \AAA$ satisfying the conclusion in the just referenced lemma.\smallskip
 	
Take now any $x\in F_p$. The fact $x\in F_p$ implies the existence of $x_0\in \AAA_0^\dd(p)$ such that $x=p+x_0$. Therefore $\mathbf{1}-x=\mathbf{1}-p-x_0$. Since both $\mathbf{1}-p$ and $x_0$ lie in $\AAA_0^\dd(p)$, we get $\mathbf{1}-x\in \AAA_0^\dd(p) \cap \AAA=L$.\smallskip

For $\lambda\in\Lambda$ define $\Theta_\lambda(x)\df \mathbf{1}-e_\lambda+U_{d_\lambda}(x)$ ($x\in F_p$). Notice that $e_\lambda=U_{d_\lambda}(\mathbf{1})$, so $\Theta_\lambda(x)=\mathbf{1}-U_{d_\lambda}(\mathbf{1}-x)$. Since $\mathbf{1}-x\in L$, we get $U_{d_\lambda}(\mathbf{1}-x)\to \mathbf{1}-x$ in norm (cf. Lemma \ref{L_Existence_of_suitable_approximate_units}), so $\Theta_\lambda(x)\to x$ in norm.\smallskip

To check additional properties of $\Theta_{\lambda}$, let us observe that for $x,y\in F_p$ we have  $$\left\|\Theta_\lambda(x)-\Theta_\lambda(y)\right\|=\left\|U_{d_\lambda}(x)-U_{d_\lambda}(y)\right\|=\left\|U_{d_\lambda}(x-y)\right\|\le \|d_\lambda\|^2\|x-y\|\le\|x-y\|,$$ so $\Theta_\lambda$ is contractive. Since $d_\lambda\ge0$ and $\|d_\lambda\|,\|-x\|\le 1$, we get from the identity  $\Theta_\lambda (x) =\mathbf{1}-U_{d_\lambda}(\mathbf{1}-x)$ ($x\in F_p$) and \cite[Lemma 3.4]{BeCuFerPe2018} that $\|\Theta_\lambda(x)\|\le 1$, so $\Theta_\lambda$ maps $F_p$ into the closed unit ball of $\AAA$.\smallskip

Fix $x\in F_p$. Since $d_\lambda, \jed-x\in L\subset \AAA_0^\dd(p)$, we deduce, via Peirce arithmetic, that $U_{d_\lambda}(\mathbf{1}-x)\in \AAA_0^\dd(p)\cap \AAA = L$. Therefore $$\Theta_\lambda (x) = \mathbf{1}-U_{d_\lambda}(\mathbf{1}-x)=p+(\mathbf{1}-p)-U_{d_\lambda}(\mathbf{1}-x)\in p+\AAA_0^\dd(p).$$ Combining the latter conclusion with $\|\Theta_\lambda(x)\|\le 1$ and $\Theta_\lambda(x)\in \AAA$ we get $$\Theta_\lambda(x)\in \AAA\cap \B_{\AAA^\dd}\cap (p+\AAA_0^\dd(p))=\AAA\cap (p+\B_{\AAA_0^\dd(p)})=F_p$$ for all $x\in F_p$.\smallskip

To simplify the notation, and in order to find the convex set $K_{\lambda}$ in the statement of the proposition, we fix a $\lambda\in\Lambda$ and we denote $e\df e_\lambda$, $f\df f_\lambda$, $\Theta\df\Theta_\lambda$.\smallskip

Set $s\df \mathbf{1}-f$. Since $0\le f\le \mathbf{1}$, $0\le s\le \mathbf{1}$. Since $f\in \AAA_0^\dd(p)$, $s=\mathbf{1}-f = p + (\mathbf{1}-p)-f\in p+\AAA_0^\dd(p)$. Therefore $s\in F_p$. Let $r\df r(s)$ be the range projection of $s$. Define $F\df\{x\in S_\AAA:r\circ x=r\},$ which is clearly non-empty.\smallskip

In order to make the following part more reader-friendly, we will extract the following statement as a claim:\smallskip

\noindent\textbf{\emph{Claim}}: {\it $F$ is a norm-closed face of $S_\AAA$, $$F=(r+\B_{\AAA_0^\dd(r)})\cap \AAA=\{x\in S_\AAA:P_2(r)(x)=r\}=\{x\in S_\AAA:x\circ s=s\},$$ and $\Theta(F_p)\subset F\subset F_p$.}

\begin{inproof}
	Namely, since $r$ is a projection, $r\circ x=\{r,r,x\}$, so $r\circ x=r$ iff $P_2(r)(x)=r$ and $P_1(r)(x)=0$. It immediately follows from this that $F=(r+\B_{\AAA_0^\dd(r)})\cap \AAA$. We also note that from \cite[Lemma 1.6]{FriRu85} if $x\in S_{\AAA^\dd}$ with $P_2(r)(x)=r$, then $P_1(r)(x)=0$. Therefore $x\in F$ if and only if $x\in S_\AAA$ and $P_2(r)(x)=r$.\smallskip
	
	It is easy to see that $F$ is closed and convex. Now let $x,y\in S_\AAA$ and $t\in(0,1)$ such that $tx+(1-t)y\in F$. Let $x_2\df P_2(r)(x)$ and  $y_2\df P_2(r)(y)$. By applying $P_2(r)$ we have $tx_2+(1-t)y_2=r$. Since Peirce projections are contractive, $x_2,y_2\in\B_{\AAA_2^\dd(r)}$. Having in mind that $r$ is the unit element in $\AAA_2^\dd(r)$, and hence it is an extreme point of the closed unit ball of this space, we deduce that $x_2=r=y_2$. Therefore $x,y\in F$ and we showed that $F$ is a face.\smallskip
	
	Let us prove the final equivalent definition of $F$. For this purpose, we observe that if $x\in F$, then $x=r+x_0$, where $x_0 = P_0(r) (x)\in \AAA_0^\dd(r)$. By combining the fact $s\in \AAA_2^\dd(r)$ with Peirce arithmetic, we arrive at $$s\circ x=\{s,\mathbf{1},x\}=\{s,r+(\mathbf{1}-r),r+x_0\}=\{s,r,r\}=s.$$ Consequently, $F\subset \{x\in S_\AAA:x\circ s=s\}$. For the reciprocal inclusion take $x\in \AAA$ such that $x\circ s=s$. Since $s$ is self-adjoint, we also have $x^\ast\circ s=s$. Put $h\df\frac{x+x^\ast}{2}$, $k\df\frac{x-x^\ast}{2}$. Then $h$ is self-adjoint, $k$ is skew-symmetric, $x=h+k$, $s\circ h=s$ and $s\circ k=0$. Having in mind that $s$ is positive we are in a position to apply \cite[Lemma 4.1]{BurFerGarPe09} to deduce that $s\perp k$ and hence $r = r(s)\perp k$ (cf. \eqref{eq reformulations of orthogonality}). On the other hand, since $(h-\mathbf{1})\circ s=0$ we can apply \cite[Lemma 4.1]{BurFerGarPe09} again to conclude that $h-\mathbf{1}\perp s$ and $h-\mathbf{1}\perp r$. Therefore $x=h+k=r+(\mathbf{1}-r)+(h-\mathbf{1})+k\in r+\AAA_0^\dd(r)$. Therefore $F= \{x\in S_\AAA:x\circ s=s\}$.\smallskip
	
	Note that since $f\circ d=d$, $f\circ e=e$ and $s=\mathbf{1}-f$, we get $s\circ d=s\circ e=0$. Therefore $\{d,d,s\}=0,$ or equivalently $d\perp s$ or $d\perp r$. Therefore, from Peirce arithmetic $U_d(x)\perp r$ for all $x\in \AAA$. Since $e\circ s=0$, we have $\mathbf{1}-e\in F$, so for any $x\in F_d$ we get $\Theta(x)=\mathbf{1}-e+U_d(x)\in r+\AAA_0^\dd(r)$, which combined with the fact $\|\Theta(x)\|\le 1$ gives $\Theta(F_d)\subset F$.\smallskip
	
	By construction $s\in F_p$, and hence there is $s_0\in \AAA_0^\dd(p)$ satisfying $s=p+s_0$. It is not hard to see that $p\perp s_0$ implies that $r=r(s)=r(p)+r(s_0)=p+r(s_0)$. Take now any $x\in F$ and write $x=r+x_0$, where $x_0 = P_0(r) (x)\perp r$. This also implies that $x_0\perp p$ (recall that $r\geq p$), and consequently $x=r+x_0=p+r(s_0)+x_0\in F_p$.\smallskip
\end{inproof}

If $r\in \AAA$, we can set ${K_\lambda=}K\df F=r+\B_{\AAA_0(r)} = r_{\lambda}+\B_{\AAA_0(r_{\lambda})}$. In this case, $r$ being a projection in $\AAA$ implies that $\mathcal{R}_{\lambda}:=\AAA_0(r)=\AAA_2(\mathbf{1}-r)$ is a unital JB$^*$-algebra, and $K$ is a proper norm closed face of $\mathcal{B}_{\AAA}$. Therefore  $\B_{\AAA_0(r)}$ (and consequently $K$) satisfies the strong Mankiewicz property (cf. \eqref{eq every convex body in a unital JB$^*$-algebra satisfies the strong Mankiewicz property}). Since $\Delta(K)$ is an intersection face, and hence convex (cf. Lemma \ref{l intersection faces MoriOzawa L8 and semi-exposed faces FerGarPeVill}), the restriction $\Delta\zuz_{K}: K\to \Delta(K)$ is affine.\smallskip

The most probable case is that $r\notin \AAA,$ a fact we assume is true from now on. Define $$C\df\{a\in \AAA:a\circ r=\gamma r\text{ for some }\gamma\in\CC\}=\left(\CC r\oplus \AAA_0^\dd(r)\right)\cap \AAA$$ (we recall that $r\circ a = \{r,r,a\} = P_2(r) (a) + \frac12 P_1(r) (a)$).\smallskip

Our next goal is to show that \begin{equation}\label{eq contractivity of Jordan product by 1-r} \hbox{$\|r\circ x\|\le\|(\mathbf{1}-r)\circ x\|,$ for any $x\in C$.}
\end{equation} Let $x=\gamma r+x_0$, where $x_0\in A_0^\dd(r)$. If $\gamma=0$, the inequality if clearly true. Otherwise we can assume $\gamma=1$. We want to show that $1\le\|x_0\|$. Arguing by contradiction, we assume that $\|x_0\|<1$. Since $x=r+x_0$, $r\perp x_0$ and $r$ is a projection, $x^n=r+x_0^n$ for all $n\in\N$. But then $\|r-x^n\|=\|x_0^n\|\le\|x_0\|^n$. From $\|x_0\|<1$ we get $\|x_0\|^n\to0$, so $x^n\to r$ in norm. But $x^n\in \AAA$ and $\AAA$ is closed in norm in $\AAA^\dd$, which contradicts that $r\notin \AAA$.\smallskip
	
Notice that for any $x_0\in \AAA_0^\dd(r)$ there is at most one $\gamma\in\CC$ such that $\gamma r+x_0\in \AAA$, otherwise we would get $r\in \AAA$. Therefore there exists a well defined mapping $\psi:(\mathbf{1}-r)\circ C\to \CC$ given by  $\psi((\mathbf{1}-r)\circ x)r=r\circ x$ for any $x\in C$. The mapping $\psi$ is obviously linear, and due to the inequality $\|r\circ x\|\le\|(\mathbf{1}-r)\circ x\|$ we have $\|\psi\|\le 1$. (Actually $\|\psi\|=1$, since $\mathbf{1}\in C$.) Also if $x,y\in C$ with $x=\gamma_x r+x_0$, $y=\gamma_y r+y_0$, where $x_0 = P_0(r) (x),y_0 = P_0(r)(y)\perp r$, then $x\circ y=\gamma_1\gamma_2 r + x_0\circ y_0$. This proves that $\psi(x_0\circ y_0)=\gamma_1\gamma_2=\psi(x_0)\psi(y_0)$ (i.e., $\psi$ is Jordan multiplicative). Moreover, if $x=\gamma r+x_0$ with $x_0\perp r$, then $x^\ast=\overline{\gamma}r+x_0^\ast$, so $\psi(x_0^\ast)=\overline{\psi(x_0)}$ (that is, $\psi$ is a Jordan $^*$-homomorphism). Actually, the mapping $C\to (\mathbf{1}-r)\circ C$, $x\mapsto  (\mathbf{1}-r)\circ x = P_0(r) (x)$ is a Jordan $^*$-homomorphism. All together implies that $\mathcal{R}_{\lambda}:=C_\R^\psi\df\psi^{-1}\left(\R\right)$ is a unital real JB$^*$-algebra, and by Proposition~\ref{p strong Mankiewicz for the real JB$^*$-algebras associated with a Jordan multiplicative functional}, its closed unit ball, $\B_{C_\R^\psi},$ satisfies the strong Mankiewicz property.\smallskip

Our next goal will consist in proving that the Jordan multiplication operators by $\mathbf{1}-r$ and by $\mathbf{1}-s$ are norm-preserving mappings when restricted to $C$. Namely, take $x\in C$ with $x=\gamma r+x_0$, where $x_0= (\jed-r)\circ x= P_0(r) (x)\perp r$. Since $r=r(s^\frac{1}{2})\geq s$, $x_0\perp s^\frac{1}{2}$, and hence $s\circ x_0=\{s^\frac{1}{2},s^\frac{1}{2},x_0\}=0$, we get $x\circ (\jed-s)= \gamma(r-s)+x_0$. Since $\|r-s\|\le\|r\| =1,$ $|\gamma|\le\|x_0\|$ (cf. \eqref{eq contractivity of Jordan product by 1-r}), we derive that \begin{equation}\label{eq isometric property of 1-s Jordan} \|x\|=\max\{\|\gamma r\|,\|x_0\|\}=\|x_0\|=\max\{\|\gamma (r-s)\|,\|x_0\|\}=\|x\circ(\jed-s)\|,
\end{equation} consequently, $\Phi_0(s): x\mapsto  x\circ(\jed-s)$ is an isometry from $C$ onto $(\jed-s)\circ C$.\smallskip
	
The set $$K_\lambda=K\df s+(\jed-s)\circ \B_{C_\R^\psi} = s_{\lambda}+(\jed-s_{\lambda})\circ \B_{C_\R^\psi},$$ which is clearly a closed convex set. By applying that $\Phi_0(s)$ and the translation by the element $s$ both are isometries, we deduce that $K_\lambda$ is truly isometrically isomorphic, via a surjective affine isometry, to the closed ball of the unital real JB$^*$-algebra $\mathcal{R}_{\lambda} ={C_\R^\psi},$ which satisfies the strong Mankiwicz property (cf. Proposition~\ref{p strong Mankiewicz for the real JB$^*$-algebras associated with a Jordan multiplicative functional}). Therefore, in order to finish the proof it suffices to show that $\Theta(F_p)\subset K\subset F_p$ and that $ \Delta(K)$ is affine.\smallskip

We shall split the arguments in a couple of auxiliary claims.\smallskip

\noindent\textbf{\emph{Claim 1}}\begin{equation}\label{claim 1 in thm 20}
	K=S_\AAA\cap\{\gamma r+(1-\gamma)s+x_0:\gamma\in[-1,1], x_0\in\AAA_0^\dd(r)\}.
\end{equation}

\begin{inproof}
	If $x=\gamma r+(1-\gamma)s+x_0$ with $x\in S_\AAA$, $\gamma\in[-1,1]$ and $x_0\in\AAA_0^\dd(r)$, then $x=s+(\jed-s)\circ(\gamma r+x_0)$. Since $\gamma r+x_0=x-(1-\gamma)s\in \AAA$, we easily see that $\gamma r+x_0\in C_\R^\psi$. Having in mind that $x_0\perp \gamma r+(1-\gamma)s$ we have $\|x_0\|\le \|x\|=1$. Then $\|\gamma r+x_0\|=\max\{|\gamma|,\|x_0\|\}\le 1$. From this we get $x\in K$.\smallskip

	Take now $x\in K$. Then there is $y=\gamma r+x_0 \in\B_\AAA$ with $\gamma\in\R$ and  $x_0\perp r$ such that $x=s+(\jed-s)\circ y$. This immediately implies that $x\in \AAA$ and $x=\gamma r+(1-\gamma)\circ s+x_0$. 
	Since $|\gamma|\le\|y\|\le1$, we get $\gamma\in[-1,1]$. From $0\le s\le\mathbf{1}$ and $r=r(s)$ we get $0\le s\le r$, so $-r\le 2s-r\le r$, which implies $\|2s-r\|\le\|r\|=1$, so $\gamma r+(1-\gamma)s=\frac{\gamma+1}{2}r+\left(1-\frac{\gamma+1}{2}\right)(2s-r)$ has norm at most one for all $\gamma\in[-1,1]$. On the other hand, as we showed before, $s=p+s_0$ and $r= r(s)=p+r(s_0)$ for some $s_0\perp p$. This implies that $\gamma r+(1-\gamma)s=p+\gamma r(s_0)+(1-\gamma)s_0\in p+\AAA_0^\dd(p)$. Therefore $\|\gamma r+(1-\gamma)s\|\ge\|p\|=1,$ which gives  $\|\gamma r+(1-\gamma)s\|=1$ for all $\gamma\in[-1,1]$. Finally, $\|x\|= \max\{\|\gamma r+(1-\gamma)s\|, \|x_0\|\}=1$,
	because $\|x_0\|\le\|y\|=1$. 
\end{inproof}

Now, if in the identity \eqref{claim 1 in thm 20} we take $\gamma=1$ we get $F\subset K$. Since $\Theta(F_p)\subset F$, we have $\Theta(F_p)\subset K$. We know that given $s=p+s_0$, $r=p+r_0$ in $F_p$, where $s_0 = P_0(p)(s),r_0 = P_0(p)(r)\in\AAA_0^\dd(p)$ (specifically $r_0=r(s_0)$), then for any $x_0\in\AAA_0^\dd(r)$ we also have $x_0\in\AAA_0^\dd(p)$, so $\gamma r+(1-\gamma)s+x_0=p+\gamma r_0+(1-\gamma)s_0+x_0,$. That gives us $K\subset F_p$.\smallskip

All what remains is to show that $\Delta(K)$ is convex. For $\gamma\in[-1,1]$, define $h_\gamma(t)\df t+\gamma(1-t)$. Notice that for any $\alpha\in[0,1]$ and any $\gamma_1,\gamma_2\in[-1,1]$ one has $\alpha h_{\gamma_1}+(1-\alpha)h_{\gamma_2}=h_{\alpha\gamma_1+(1-\alpha)\gamma_2}$. For any set $S\subset \AAA$ denote by $\mathcal{N}_\delta(S)$ the $\delta$-neighbourhood of $S$ in $\AAA$. Denote $r_{k,m}^i\df \chi_{\left[\frac{2k-2+i}{2m},\frac{2k-1+i}{2m}\right]}(s)$, where $m\in\N$, $i=1,2$, $k=1,\ldots,m$, and set $$G_m^i(\gamma)\df F_p\bigcap\left(\bigcap_kF\left(r_{k,m}^i,h_\gamma\left(\frac{k}{m}\right)\right)\right)$$ and $$H_m^i(\gamma)\df F_p\cap\mathcal{N}_{\frac{1}{m}}\left(G_m^i(\gamma)\right),$$ where the intersection is over those $k=1,2,\dots,m$  for which $r_{k,m}^i\ne0$. Clearly $r_{k,m}^i\le r$ for all $i,m,k$ as above. Theorem~\ref{t lemma 17 for triples}, combined with the identity $\alpha h_{\gamma_1}+(1-\alpha)h_{\gamma_2}=h_{\alpha\gamma_1+(1-\alpha)\gamma_2}$ and \cite[Lemma 10]{MoriOza2018} leads to $$\alpha G_m^i(\gamma_1)+(1-\alpha)G_m^i(\gamma_2)\subset G_m^i(\gamma_3)$$ and consequently $$\alpha H_m^i(\gamma_1)+(1-\alpha)H_m^i(\gamma_2)\subset H_m^i(\gamma_3)$$ for any $\gamma_1,\gamma_2\in[-1,1]$, $\alpha\in[0,1]$, and $\gamma_3\df\alpha\gamma_1+(1-\alpha)\gamma_2$.\smallskip

In order to have a more friendly proof, we insert next the following:\smallskip

\noindent\textbf{\emph{Claim 2}} Denoting $K(\gamma)\df\{x\in F_p:x\circ r=h_\gamma(s)\}$ we have
\begin{equation}\label{Claim 2 in thm 20}
	 K(\gamma)=\bigcap_{m\in\N}\left(H_m^1(\gamma)\cap H_m^2(\gamma)\right).
\end{equation}

\begin{inproof}($\subset$) Take any $x\in K(\gamma)$. Let us observe that $x\circ r =\{r,r,x\}= P_2(r) (x) + \frac12 P_1(r) (x)$, and hence $r\circ x = h_{\gamma} (s)\in \AAA^{**}_2 (r)$ implies that $x = h_{\gamma} (s) + P_0 (r) (x)$.  Define $g_m$ to be the continuous function defined as $$g_m(t):=\left\{ \begin{array}{ll}
			\gamma , & \hbox{if $t = 0$,} \\
			h_\gamma\left(\frac{k}{m}\right), & \hbox{if $t\in\left[\frac{2k-1}{2m},\frac{2k}{2m}\right]$,} \\
			\hbox{affine }, & \hbox{ if $t\in\left[\frac{2k-2}{2m},\frac{2k-1}{2m}\right]$}. 
		\end{array} \right.$$ Then since $h_\gamma\left(\frac{k}{m}\right)-h_\gamma\left(\frac{2k-1}{2m}\right)=\frac{1-\gamma}{2m}$, we have $\|g_m-h_\gamma\|_\infty\le\frac{1}{m}$. By definition we also have $(g_m-h_\gamma)(0)=0$, and thus $(g_m-h_\gamma)(s)\in \AAA$ and $\|(g_m-h_\gamma)(s)\|\le\frac{1}{m}$. Denote $y\df x+(g_m-h_\gamma)(s)$. To show $x\in H_m^1(\gamma)$ it suffices to prove that $y\in G_m^1(\gamma)$.\smallskip
		
		We know that $y\in \AAA$.  By observing that $y = (\jed-r) \circ x + g_m(s) = P_0(r) (x) + g_m(s)$ we deduce that $\|y\| = \max\{ \|P_0(r) (x)\|, \|g_m(s)\|\} = 1$. Again, $s=p+s_0$ for some $s_0\perp p$. Since $p$ is a projection, this implies that $P(s)=P(1) p + P(s_0)$ for any polynomial $P$ with zero constant term, and by the Stone--Weierstrass theorem $f(s)=f(1) p + f(s_0)$ for any $f\in C([0,1])$ with $f(0)=0$. Since $(g_m-h_\gamma)(1)=0$, we get $(g_m-h_\gamma)(s)=(g_m-h_\gamma)(s_0)\in\AAA_0^\dd(p)$. Since $x\in F_p$, we get also $y\in F_p$.\smallskip

Now, the equality $x\circ r=h_\gamma(s)$ implies that  $y\circ r=g_m(s)$. Moreover, the identities $$g_m\cdot \chi_{\left[\frac{2k-1}{2m},\frac{2k}{2m}\right]}=h_\gamma\left(\frac{k}{m}\right)\cdot \chi_{\left[\frac{2k-1}{2m},\frac{2k}{2m}\right]},$$ and $y\circ r=g_m(s)$ give $y\circ r_{k,m}^1=h_\gamma\left(\frac{k}{m}\right)r_{k,m}^1$, or equivalently, $y\in F\left(r_{k,m}^1,h_\gamma\left(\frac{k}{m}\right)\right)$ for all $k$. Therefore $y\in G_m^1(\gamma)$ and $x\in H_m^1(\gamma)$. Summarizing $K(\gamma)\subset H_m^1(\gamma)$. Analogously, $K(\gamma)\subset H_m^2(\gamma)$.\smallskip
				
($\supset$) Take $x\in\bigcap_{m\in\N}\left(H_m^1(\gamma)\cap H_m^2(\gamma)\right)$. Then for any $i,m$ there exists $y^i_m\in G_m^i(\gamma)$ such that $\|x-y_m^i\|\le \frac{1}{m}$. Since $$y_m^i\in F\left(r_{k,m}^i,h_\gamma\left(\frac{k}{m}\right)\right)=h_\gamma\left(\frac{k}{m}\right)r_{k,m}^i+\B_{\AAA_0^\dd(r_{k,m}^i)},$$ $y_m^i\circ r_{k,m}^i=h_\gamma\left(\frac{k}{m}\right)r_{k,m}^i $. Write $r_m^i\df\sum_{k=1}^mr_{k,m}^i$. Then $$y_m^i\circ r_m^i=\sum_{k=1}^my_m^i\circ r_{k,m}^i=\sum_{k=1}^mh_\gamma\left(\frac{k}{m}\right)r_{k,m}^i=\left(\sum_{k=1}^mh_\gamma\left(\frac{k}{m}\right)\chi_{\left[\frac{2k-2+i}{2m},\frac{2k-1+i}{2m}\right]}\right)(s).$$ As before, $$\left\|\sum_{k=1}^mh_\gamma\left(\frac{k}{m}\right)\chi_{\left[\frac{2k-2+i}{2m},\frac{2k-1+i}{2m}\right]}-\sum_{k=1}^mh_\gamma\cdot\chi_{\left[\frac{2k-2+i}{2m},\frac{2k-1+i}{2m}\right]}\right\|_\infty\le\frac{1}{m},$$ which implies that  $$\|h_\gamma(s)\circ r_m^i-y_m^i\circ r_m^i\|\le\frac{1}{m}.$$ Therefore from $\|x-y_m^i\|\le \frac{1}{m}$ we deduce that  $\left\|(h_\gamma(s)-x)\circ r_m^i\right\|\le\frac{2}{m}$ for $i=1,2$.\smallskip
		
		Let $r_m\df r_m^1\vee r_m^2$. Since $r_m^i=\sum_{k=1}^m\chi_{\left[\frac{2k-2+i}{2m},\frac{2k-1+i}{2m}\right]}(s)$, we can easily see that $r_m=\chi_{\left[\frac{1}{m},1\right](s)}$ and $r_m=r_m^1+r_m^2-r_m^1\circ r_m^2=r_m^2+r_m^1\circ (\mathbf{1}-r_m^2)$. Notice that $r_m^1\circ (\jed-r_m^2)=\chi_{\left[\frac{1}{2m},\frac{2}{2m}\right)}+\sum_{k=2}^m\chi_{\left(\frac{2k-2+i}{2m},\frac{2k-1+i}{2m}\right)}$ is a projection and $0\le r_m^1\circ (\jed-r_m^2)\le r_m^1$ (all these projections lie in a JBW$^*$-subalgebra of $\AAA^{**}$ which is isometrically isomorphic to a commutative von Neumann algebra).\smallskip 
		
		Since $r_m^1$ and $r_m^1\circ (\jed-r_m^2)$ are projections with $0\le r_m^1\circ (\jed-r_m^2)\le r_m^1$ and $\left\|(h_\gamma(s)-x)\circ r_m^1\right\|\le\frac{2}{m}$, Lemma~\ref{L_Jordan_product_estimate_for_general_elements} implies that  $\|(h_\gamma(s)-x)\circ \left(r_m^1\circ (\mathbf{1}-r_m^2)\right)\|\le \frac{12}{m}$. Therefore $$\left\|(h_\gamma(s)-x)\circ r_m\right\|\le \left\|(h_\gamma(s)-x)\circ \left(r_m^1\circ (\mathbf{1}-r_m^2)\right)\right\|+\left\|(h_\gamma(s)-x)\circ r_m^2\right\|\le\frac{14}{m}.$$
		
		Since $s\ge0$, $r=r(s)=\chi_{(0,1]}(s)$, and $r_m\nearrow r$ in the strong$^*$ topology of $\AAA^{**}$. Therefore $x\circ r_m\to x\circ r$ in the strong$^*$ topology of $\AAA^{**}$. \smallskip
		
		On the other hand, it is known that for each $z\in \AAA^{**}$ we have  $\|z\|=\sup_{\phi\in S_{\AAA^\ast}}\|z\|_\phi$ (see \cite[Proposition 1.2]{BarFri87Groth} or \cite[Proposition 5.10.60]{CabGarPalVol2}). Moreover, since the strong$^*$ topology of  $\AAA^{**}$ is the topology generated by all the preHilbertian seminorms $\|\cdot\|_{\phi}$ with $\phi\in  S_{\AAA^\ast}$ (cf. \cite[Definition 3.1]{BarFri90}), $$\left\|(h_\gamma(s)-x)\circ r_m\right\|_\phi\le\left\|(h_\gamma(s)-x)\circ r_m\right\|\le\frac{14}{m}$$ for any $\phi\in S_{\AAA^\ast}$ (see \cite[Proposition 1.2]{BarFri87Groth}), and $(h_\gamma(s)-x)\circ r_m\to (h_\gamma(s)-x)\circ r$ in the strong$^*$ topology, and hence $\|(h_\gamma(s)-x)\circ r_m\|_\phi\to \|(h_\gamma(s)-x)\circ r\|_\phi$ for any $\phi\in S_{\AAA^\ast}$, we deduce that $\|(h_\gamma(s)-x)\circ r\|_\phi=0$ for all $\phi\in S_{\AAA^\ast}$, and thus $$\|(h_\gamma(s)-x)\circ r\|=\sup_{\phi\in S_{\AAA^\ast}}\|(h_\gamma(s)-x)\circ r\|_\phi=0,$$ equivalently, $x\circ r=h_\gamma(s)\circ r=h_\gamma(s)$. This means that $x\in K(\gamma)$, which is what we wanted. 
\end{inproof}

Let us now conclude the proof of the proposition.  Since $h_\gamma(s)=s+\gamma(r-s)=\gamma r+(1-\gamma)s$, an element $x\in S_\AAA$ satisfies $x=\gamma r+(1-\gamma)s+x_0$, where $x_0\perp r$, if and only if $x\circ r=h_\gamma(s)$. From this observation combined with Claim 1 \eqref{claim 1 in thm 20} we get \begin{equation}\label{identity K as union of Kgamma} K=\bigcup_{\gamma\in[-1,1]}K(\gamma).
\end{equation}

Let us focus on the sets $G_m^i(\gamma)$ and $H_m^i(\gamma).$ We know that $F_p$ is a norm closed proper face of $S_{\AAA}$ and $\Delta(F_p)$ is an intersection face in $S_{Y}$ (cf. Lemma \ref{l intersection faces MoriOzawa L8 and semi-exposed faces FerGarPeVill}). By Theorem \ref{t lemma 17 for triples} $$F\left(r_{k,m}^i,h_\gamma\left(\frac{k}{m}\right)\right) = \left(F_{r_{k,m}^i}\right)_{h_\gamma\left(\frac{k}{m}\right)} ,$$ which combined with Lemma \ref{l 11 MoriOzawa}$(b)$ assures that $$\Delta\left(F\left(r_{k,m}^i,h_\gamma\left(\frac{k}{m}\right)\right)\right) = \Delta\left(\left(F_{r_{k,m}^i}\right)_{h_\gamma\left(\frac{k}{m}\right)}\right) =  \Delta\left(F_{r_{k,m}^i}\right)_{h_\gamma\left(\frac{k}{m}\right)}. $$ Therefore 
$$\begin{aligned}
\Delta\left(G_m^i(\gamma)\right) &=  \Delta\left(F_p\right) \bigcap \left(\bigcap_k \Delta\left(\left(F_{r_{k,m}^i}\right)_{h_\gamma\left(\frac{k}{m}\right)}\right)  \right) \\
&=  \Delta\left(F_p\right) \bigcap \left(\bigcap_k \Delta\left(F_{r_{k,m}^i}\right)_{h_\gamma\left(\frac{k}{m}\right)}  \right),
\end{aligned}$$ and $$\Delta\left(H_m^i(\gamma)\right) =  \Delta\left(F_p\right) \cap\mathcal{N}_{\frac{1}{m}}\left(\Delta\left(G_m^i(\gamma)\right)\right).$$ 

Given $\alpha\in [0,1]$, Lemma \ref{l 11 MoriOzawa}$(c)$ (see also \cite[Lemma 11]{MoriOza2018}) assures that 
$$\begin{aligned}
&\alpha \Delta\left(G_m^i(\gamma_1)\right) + (1-\alpha) \Delta\left(G_m^i(\gamma_2)\right) \\ &\subseteq \Delta\left(F_p\right) \bigcap \left(\bigcap_k \alpha \Delta\left(\left(F_{r_{k,m}^i}\right)_{h_{\gamma_1}\left(\frac{k}{m}\right)}\right) +(1-\alpha) \Delta\left(\left(F_{r_{k,m}^i}\right)_{h_{\gamma_2}\left(\frac{k}{m}\right)}\right) \right) \\
&\subseteq \Delta\left(F_p\right) \bigcap \left(\bigcap_k  \Delta\left(\left(F_{r_{k,m}^i}\right)_{ \alpha  h_{\gamma_1}\left(\frac{k}{m}\right) +(1-\alpha) h_{\gamma_2}\left(\frac{k}{m}\right) }\right)  \right) \\
&= \Delta\left(F_p\right) \bigcap \left(\bigcap_k  \Delta\left(\left(F_{r_{k,m}^i}\right)_{ h_{ \alpha \gamma_1 +(1-\alpha) \gamma_2 } \left(\frac{k}{m}\right) }\right)  \right) \\
&= \Delta\left(F_p\right) \bigcap \left(\bigcap_k  \Delta\left(\left(F_{r_{k,m}^i}\right)_{ h_{\gamma_3 } \left(\frac{k}{m}\right) }\right)  \right) = \Delta\left(G_m^i(\gamma_3)\right)\\
\end{aligned}$$ with $\gamma_3 = \alpha \gamma_1 +(1-\alpha) \gamma_2$.\smallskip

We consequently have 
\begin{equation}\label{eq convex comb Delta Hi}
\begin{aligned}
	&\alpha \Delta(H_m^i(\gamma_1))  + (1-\alpha) \Delta(H_m^i(\gamma_2)) \\ 
	&\subseteq \Delta\left(F_p\right) \bigcap \Big(\alpha \mathcal{N}_{\frac{1}{m}}\left(\Delta\left(G_m^i(\gamma_1)\right)\right) + (1-\alpha) \mathcal{N}_{\frac{1}{m}}\left(\Delta\left(G_m^i(\gamma_2)\right)\right)\Big)\\ 
	&\subseteq \Delta\left(F_p\right) \bigcap \mathcal{N}_{\frac{1}{m}}\Big( \alpha \Delta\left(G_m^i(\gamma_1)\right) + (1-\alpha) \Delta\left(G_m^i(\gamma_2)\right)\Big) \\ 
	&\subseteq \Delta\left(F_p\right) \bigcap \mathcal{N}_{\frac{1}{m}}\Big(\Delta\left(G_m^i(\gamma_3)\right) \Big) = \Delta(H_m^i(\gamma_3))
\end{aligned}
\end{equation} with $\gamma_3 = \alpha \gamma_1 +(1-\alpha) \gamma_2$.\smallskip

By Claim 2 \eqref{Claim 2 in thm 20}, $$
K(\gamma)=\bigcap_{m\in\N}\left(H_m^1(\gamma)\cap H_m^2(\gamma)\right).$$ So, given $\gamma_1,\gamma_2\in [-1,1]$ and $\alpha\in [0,1]$ we deduce from \eqref{eq convex comb Delta Hi} that $$\begin{aligned}
& \alpha \Delta(K(\gamma_1))+(1-\alpha)\Delta(K(\gamma_2))\\
&\subseteq \bigcap_{m\in\N} \Big( \alpha \Delta\Big(H_m^1(\gamma_1)\cap H_m^2(\gamma_1)\Big) + (1-\alpha) \Delta\Big(H_m^1(\gamma_2)\cap H_m^2(\gamma_2)\Big) \Big)\\
&\subseteq \bigcap_{m\in\N}  \Delta\Big(H_m^1(\gamma_3)\cap H_m^2(\gamma_3)\Big) =K(\gamma_3).
\end{aligned}$$ We therefore conclude from \eqref{identity K as union of Kgamma} that $\Delta(K)$ is convex, which finishes the proof.
\end{proof}

A strengthened version of the previous proposition will be obtained in the next section.  

\section{Minimal tripotents in the second dual of a unital JB$^*$-algebra are always dominated by unitaries in the principal component}\label{sec: domination of minimal tripotents in the second dual by unitaries in the main component}

The bidual, $E^{**},$ of each JB$^*$-triple, $E$, admits an \emph{atomic decomposition} in the form $E^{**} = \mathcal{A}\oplus^{\infty} \mathcal{N}$, where $\mathcal{A}$ and $\mathcal{N}$ are two orthogonal weak$^*$-closed ideals of $E^{**}$ --called the \emph{atomic} and \emph{non-atomic} part of $E^{**},$ respectively--, $\mathcal{A}$ coincides with the weak$^*$-closed linear span of all minimal tripotents in $E^{**}$ and its predual coincides with the norm closure of the linear span of all extreme points of the closed unit ball of $E^{*}$, the closed unit ball of the predual of $\mathcal{N}$ has no extreme points and $\mathcal{N}$ contains no minimal tripotents \cite[Theorems 1 and 2]{FriRu85}. It is further known that $\mathcal{A}$ is isometrically isomorphic, as JB$^*$-triple, with an $\ell_{\infty}$-sum of a family of Cartan factors, which paves the way to embed $E$ as a JB$^*$-subtriple of an $\ell_{\infty}$-sum of Cartan factors \cite[Proposition 2 and Theorem 1]{FriRu86}.\smallskip

Let $E$ be a JB$^*$-triple. Following the terminology in JB$^*$-triple theory \cite{BunChu92}, we denote by $K_0(E)$ the Banach subspace of $E$ generated by the minimal tripotents of $E$. The space $K_0(E)$ can be described as the largest inner ideal of $E$ which is also an inner ideal of $E^{**}$ \cite[Corollary 3.5]{BunChu92}. The algebraic linear span of all minimal tripotents in $E$ is called the \emph{socle} of $E$ (soc$(E)$ in short). The socle of $E$ need not be a closed subspace. That is the case of the socle of $B(H)$, which coincides with the subspace, $\mathcal{F} (H)$, of all finite rank operators, and it is not closed when $H$ is infinite dimensional. However, if a JB$^*$-triple $E$ has finite rank (equivalently, if $E$ is a reflexive JB$^*$-triple), we have soc$(E) = E$ (cf. \cite[Proposition 4.5 and Remark 4.6]{BunChu92} or \cite{BeLoPeRo}).\smallskip

For a JB$^*$-triple $E$, the inner ideal $K_0(E^{**})$ and the socle of $E^{**}$ are better understood in terms of elementary JB$^*$-triples. Given a Cartan factor of type $j\in \{1,\ldots, 6\},$ the elementary JB$^*$-triple $K_j$ of type $j$ is defined in the following terms: $K_1 = K (H_1, H_2)$; $K_i = C \cap K(H)$ when $C$ is of type $i = 2 , 3$, where $K(H_1,H_1)$ stands for the space of all compact linear operators from $H_1$ to $H_2$ and $K(H) = K(H,H)$, and we set $K_i = C$ if the latter is of type $ 4, 5,$ or $6$. Obviously, if $K$ is an elementary JB$^*$-triple of type $j$, its bidual is precisely a Cartan factor of $j$. Let us write $\mathcal{F}_1= \mathcal{F} (H_1,H_2)$ for the non-necessarily closed ideal of $B(H_1,H_2)$ of all finite--rank linear operators from $H_1$ to $H_2$, $\mathcal{F} (H) = \mathcal{F} (H,H)$, $\mathcal{F}_j  = C\cap \mathcal{F} (H)$ if $C$ is a type $j = 2,3$ Cartan factor, and $\mathcal{F}_j =C$ in the remaining cases. Thanks to this notation, we can see that if the atomic part of a JB$^*$-triple $E$ coincides with the $\ell_{\infty}$-sum of a family of Cartan factors $\{C_j : j\in \Lambda\}$, then $K_0(E^{**})$ is precisely the ${c_0}$-sum of the corresponding elementary JB$^*$-triples $\{K_j : j\in\Lambda\},$ and the socle of $E^{**}$ is the non-necessarily closed inner ideal generated by the $\{\mathcal{F}_j : j\in \Lambda\}$. \smallskip

Lemma 2.3 and Proposition 5 in \cite{FriRu85} establish a complete description of the JB$^*$-subtriple generated by two minimal tripotents $u,v$ in a JB$^*$-triple $E$. It is shown that this subtriple is of dimension at most $4$, it is linearly generated by $u,v,P_1(u) (v)$ and $P_1(v) (u)$ and isometrically isomorphic to $\mathbb{C}$, $M_{1,2}(\mathbb{C})$, $\mathbb{C}\oplus^{\infty} \mathbb{C}$, $M_2(\mathbb{C})$ or $S_2(\mathbb{C})$, where the latter denotes the space of all symmetric complex matrices. In our first result we shall focus on the JB$^*$-subalgebra generated by a minimal tripotent in a JB$^*$-algebra.\smallskip

Let us briefly comment that the involution on each JB$^*$-algebra $\AAA$ is a conjugate-linear triple automorphism on $\AAA$. So, if $e$ is a (minimal) tripotent in $\AAA$, its adjoint $e^*$ also is a (minimal) tripotent in $\AAA$. In our next result we describe the JB$^*$-subalgebra generated by a minimal tripotent in an arbitrary JB$^*$-algebra. 

\begin{proposition}\label{p representation of the JBstar algebra generated by a minimal projection} 
	Let $\AAA$ be a JB$^*$-algebra and let $e$ be a minimal tripotent in $\AAA$. Let $\mathfrak{B}$ stand for the JB$^*$-subalgebra of $\AAA$ generated by $e$. Then $\mathfrak{B}$ is a unital JW$^*$-algebra of dimension at most $3$ and rank at most $2$, being linearly spanned by $e,e^*$ and $e\circ e^*$. Furthermore, $\mathfrak{B}$ is, in general, strictly bigger than the JB$^*$-subtriple of $\AAA$ generated by the minimal tripotents $e$ and $e^*$.
\end{proposition}

\begin{proof} To simplify the notation, let $\mathfrak{F}$ be the JB$^*$-subtriple of $\AAA$ generated by the minimal tripotents $e$ and $e^*$. Since $e,e^*\in \mathfrak{B}$ and the latter also is a JB$^*$-subtriple of $\AAA$, it follows that $\mathfrak{F} \subseteq \mathfrak{B}$. \smallskip
	
	We shall first show that $e^2$ must be a scalar multiple of $e$. Namely, since, by weak$^*$-density of $\AAA$ in $\AAA^{**}$ and the triple product of $\AAA^{**}$ is separately weak$^*$-continuous, the tripotent $e$ also is minimal in $\AAA^{**}$, we can therefore assume that $\AAA$ is a unital JB$^*$-algebra. By Peirce arithmetic $$e^2 = \{e,\mathbf{1},e\} \in \AAA_2(e) = \mathbb{C} e.$$ Therefore, there exists $\lambda\in\CC$ such that $e^2=\lambda e$. Consequently,  $(e^\ast)^2=\overline{\lambda}e^\ast$.\smallskip 
	
	Now $$e = \{e,e,e\} = 2 e\circ(e\circ e^\ast) - e^2 \circ e^*= 2 e\circ(e\circ e^\ast) - \lambda e\circ e^* ,$$ witnessing that \begin{equation}\label{eq ecirc e star circ e} e\circ(e\circ e^\ast)  =\frac{e}{2}+\frac{\lambda}{2}e\circ e^\ast, \hbox{ and } e^\ast\circ(e\circ e^\ast)=\frac{e^\ast}{2}+\frac{\overline{\lambda}}{2}e\circ e^\ast.
	\end{equation} 
	
	Next, by applying that $\lambda e = e^2 = \{e,\mathbf{1},e\}$, by the Jordan identity, we have 
	$$\begin{aligned} \overline{\lambda} e & = \{e,\{e,\mathbf{1},e\},e\} = - \{\mathbf{1}, e, \{e,e,e\}\} + 2 \{e,e, \{\mathbf{1}, e,e\}\}  \\
		&= - \{\mathbf{1}, e, e\}  + 2 \{e,e, e\circ e^*\} = -  e\circ e^*  + 2  (e\circ e^*)\circ (e\circ e^*) \\
		&+ 2 ((e\circ e^*) \circ e^*) \circ e - 2 ((e\circ e^*) \circ e) \circ e^* \\
		&= 2  (e\circ e^*)^2   -  e\circ e^* + 2\left(\frac{e^\ast}{2}+\frac{\overline{\lambda}}{2}e\circ e^\ast\right) \circ e - 2 \left(\frac{e}{2}+\frac{\lambda}{2}e\circ e^\ast\right) \circ e^* \\
		&= 2  (e\circ e^*)^2   -  e\circ e^* +\overline{\lambda} \left(e\circ e^\ast\right) \circ e - \lambda \left(e\circ e^\ast\right) \circ e^*\\
		&= 2  (e\circ e^*)^2   -  e\circ e^* +\overline{\lambda} \left(  \frac{e}{2}+\frac{\lambda}{2}e\circ e^\ast \right) - \lambda \left(\frac{e^\ast}{2}+\frac{\overline{\lambda}}{2}e\circ e^\ast\right) \\
		&= 2  (e\circ e^*)^2   -  e\circ e^* +\frac{\overline{\lambda}}{2} e - \frac{\lambda}{2} e^\ast,
	\end{aligned}$$ from which we deduce that 
	$$\begin{aligned} (e\circ e^\ast)^2 &=\frac12  e\circ e^* +\frac{\overline{\lambda}}{4} e + \frac{\lambda}{4} e^\ast.
	\end{aligned}$$
This shows that the set $\spn\{e,e^\ast,e\circ e^\ast\}$ is closed for Jordan products, and since it has dimension at most 3, it is closed in $\AAA$. Therefore the JB$^*$-subalgebra $\mathfrak{B}$ is contained in the linear span $e,e^\ast,$ and $e\circ e^\ast$, and thus it has dimension at most three. Consequently, $\mathfrak{B}$ is a finite dimensional unital JBW$^*$-algebra.\smallskip

We have already observed that $\mathfrak{F}\subseteq \mathfrak{B}$. To show that $\mathfrak{B}$ is, in general, strictly bigger than $\mathfrak{F}$, it suffices to consider the minimal tripotent $e = \frac{1}{2}  \left(\begin{matrix}	1 & 1 \\
	-1	& -1 
\end{matrix}\right)$ in $\mathfrak{A} = M_2 (\mathbb{C})$. It is easy to check that $P_1 (e) (e^*) = (\mathbf{1}-ee^*) e^* (e^*e) + (ee^*) e^* (\mathbf{1}-e^*e ) =0$ $= P_1 (e^*) (e)$, and thus the JB$^*$-subtriple of $\AAA$ generated by  $e$ and $e^*$ is two dimensional (cf. the comments before this proposition or \cite[Lemma 2.3]{FriRu85}). However $e\circ e^* =  \frac{1}{2}  \left(\begin{matrix}	1 & 0 \\
	0	& 1 
\end{matrix}\right)$, and the elements $e,e^*$ and $e\circ e^* $ are linearly independent, which implies that $\mathfrak{B}$ has dimension three.\smallskip

For the statement concerning the rank of $\mathfrak{B}$ we observe that that this rank must be finite by the finite dimensionality of $\mathfrak{B}$. Furthermore, if $\mathfrak{B}$ has rank $r\geq 3$, we can find mutually orthogonal minimal tripotents $v_1,v_2,\ldots, v_r$ in $\mathfrak{B}$. The conclusion on the dimension of $\mathfrak{B}$ proves that $r= 3,$ $\mathfrak{B}$ must linearly generated by $v_1,v_2$ and $v_3$, and by minimality of $e$, we can assume that $e = \gamma v_1$ with $|\gamma|=1$. Therefore $e, e^* = \overline{\gamma} v_1^*$ and $e\circ e^*$ can only generate a two dimensional space, which contradicts that $\mathfrak{B}$ contains three mutually orthogonal minimal tripotents.\smallskip

Finally, since $\mathfrak{B}$ is the JB$^*$-subalgebra generated by one element, the Shirshov--Cohn theorem implies that $\mathfrak{B}$ is a JW$^*$-algebra.  
\end{proof}

\begin{remark}\label{r consequences of the representation of the JB$^*$-subalgebra generated by a min trip} {\rm There are several interesting/surprising consequences of our previous proposition. As we commented before, the JB$^*$-subtriple of a JB$^*$-triple $E$ generated by two minimal tripotents can be one, two, three or four dimensional (cf. \cite[Proposition 5]{FriRu85}). However, as a consequence of the above Proposition~\ref{p representation of the JBstar algebra generated by a minimal projection}, the JB$^*$-subtriple $\mathfrak{F}$ of a JB$^*$-algebra $\AAA$ generated by a minimal tripotent $e$ and its transposed $e^*$ is at most three dimensional.  It is also worth to explore when $\mathfrak{F}$ coincides with the JB$^*$-subalgebra $\mathfrak{B}$ generated by $e$. We shall see that the following statements are equivalent:
	\begin{enumerate}[$(a)$]\item $\mathfrak{F} = \mathfrak{B}$;
		\item $e^2\neq 0$ (equivalently, $(e^*)^2 \neq 0$).
	\end{enumerate}		
	
	It follows from the arguments in the proof of Proposition~\ref{p representation of the JBstar algebra generated by a minimal projection} and \cite[Lemma 2.3]{FriRu85} that $\mathfrak{F} = \mathfrak{B}$ if and only if $e\circ e^* $ is a linear combination of $e,e^*,$ $P_1(e^*) (e)$ and $P_1(e) (e^*)$. As in the proof of Proposition~\ref{p representation of the JBstar algebra generated by a minimal projection}, $e^2 = \lambda e$ for some $\lambda\in \mathbb{C}$. Since $e$ is minimal there exists a scalar $ \mu_{_{e^*}}\in \mathbb{C}$ such that $ \mu_{_{e^*}} e = P_2(e) (e^*)$. We deduce from the Jordan identity that \begin{equation}\label{eq 1 remark 2511} \begin{aligned} \mu_{_{e^*}} e + \frac12 P_1 (e) (e^*)  &= P_2(e) (e^*) + \frac12 P_1 (e) (e^*) =  L(e,e) (e^* ) \\
			=& (e\circ e^*)\circ e^* + 	(e^*\circ e^*)\circ e - (e\circ e^*)\circ e^* \\
			&= (e^*\circ e^*)\circ e = \overline{\lambda}\  e^*\circ e.
		\end{aligned}
	\end{equation} So, if $\lambda\neq 0$, the element $e^*\circ e$ belongs to the linear span of $e$ and $P_1 (e) (e^*)$ (and, by taking adjoints, to the linear span of $e^*$ and $P_1 (e^*) (e)$). \smallskip
	
	If $\lambda = 0$ (equivalently, $e^2 =0$), it can be easily seen from \eqref{eq 1 remark 2511} that $ P_2(e) (e^*) =  P_1 (e) (e^*) = 0$. Taking adjoints we get $ P_2(e^*) (e) =  P_1 (e^*) (e) = 0$. In particular, $e^* =P_0(e) (e^*)$ is orthogonal to $e$. Under these circumstances $\mathfrak{F}$ is at most two dimensional being generated by $e$ and $e^*$. This shows that $e\circ e^*\in \mathfrak{F}$ if and only if we can write $e\circ e^*  = \alpha e +\beta e^*$ for some $\alpha,\beta\in \mathbb{C}$.  Having in mind \eqref{eq ecirc e star circ e} we deduce that $$\beta \alpha e +\beta^2 e^* =\beta e\circ e^*= e\circ(e\circ e^\ast)  = \frac{1}{2} e,$$ and similarly, $$ \alpha^2 e +\alpha \beta e^* =\alpha e\circ e^* = e^\ast\circ(e\circ e^\ast)=\frac{1}{2} e^\ast.$$ By applying that $e\perp e^*$ we arrive at $\alpha^2 = \beta^2 =0$, equivalently, $e\circ e^* =0$. However, this implies that $e =0$, which is impossible. 
}\end{remark}

In the next corollary we combine the previous conclusion with some recent results from \cite{HKP-Fin}. We recall that a tripotent $e$ in a JBW$^*$-triple $M$ is called {\em finite} if any tripotent $u\in M_2(e)$ which is complete in $M_2(e)$ is already unitary in $M_2(e)$. If every tripotent in $M$ is finite, we say that $M$ itself is {\em finite}.

\begin{corollary}\label{c representation of the JBstar algebra generated by a minimal projection in the bidual} 
Let $\AAA$ be a JB$^*$-algebra, $e$ a minimal tripotent in $\AAA^\dd$, and let $\mathfrak{B}$ be the JB$^*$-subalgebra of $ \AAA^\dd$ generated by $e$. Then
\begin{enumerate}[$(a)$]
	\item There is a $u\in \mathcal{U}(\mathfrak{B})$ such that $e\le u$;
	\item $\mathfrak{B}$ is contained in the socle of $\AAA^\dd$. In particular, the unit of $\mathfrak{B}$ is a finite rank tripotent in $\AAA^\dd$.
\end{enumerate} 
\end{corollary}

\begin{proof} $(a)$ Since, by Proposition~\ref{p representation of the JBstar algebra generated by a minimal projection}, $\mathfrak{B}$ is a finite dimensional JW$^*$-algebra, we deduce from \cite[Proposition 3.4]{HKP-Fin} that it is a finite JBW$^*$-triple. Therefore $e$ is a finite tripotent in the JBW$^*$-algebra $\mathfrak{B}$, so by \cite[Proposition 7.5]{HKP-Fin} there is a unitary element $u\in \mathfrak{B}$ such that $e\le u$.\smallskip

$(b)$ Since $e,e^*\in \hbox{soc}(\AAA^{**})$ and the latter is an algebraic inner ideal of $\AAA^{**}$, we deduce that $e\circ e^* = \{e,\mathbf{1},e^*\}\in \hbox{soc} (\AAA^\dd)$, and thus $\mathfrak{B}\subset  \hbox{soc} (\AAA^\dd)$ (cf. Proposition~\ref{p representation of the JBstar algebra generated by a minimal projection}). The rest is clear.
\end{proof}

The previous corollary together with the JB$^*$-triple version of Kadison's transivity theorem \cite{BurFerGarMarPe,Ham99} give the following geometric conclusion which is interesting by itself and might have many different applications.

\begin{theorem}\label{t boundedness  principle for minimal tripotents in the second dual by unitaries in the principal component} Let $\AAA$ be a JB$^*$-algebra and let $e$ be a minimal tripotent in $\AAA^\dd$. Then there exists a self-adjoint element $h$ in $\AAA$ satisfying $e\leq \exp(ih)$, that is, $e$ is bounded by a unitary in the principal connected component of $\mathcal{U} (\AAA)$.  
\end{theorem}   

\begin{proof} Let $\mathfrak{B}$ denote the JB$^*$-subalgebra of $\AAA^{**}$ generated by $e$. By Corollary~\ref{c representation of the JBstar algebra generated by a minimal projection in the bidual}, there exists a unitary $u\in \mathfrak{B}\subset \hbox{soc} (\AAA^{**})$ with $e\leq u$. By applying that $\mathfrak{B}$ is finite dimensional and it is contained in the socle of $\AAA^{**},$ we can find a spectral resolution of $u$ in the form $ u = \exp(i \theta_1) p_1 +\ldots+  \exp(i \theta_k) p_k$, where $p_1,\ldots, p_k$ are mutually orthogonal projections in $\mathfrak{B}$ which are minimal in $\AAA^{**}$ and $\theta_1,\ldots, \theta_k$ are real numbers. \smallskip

By 	Kadison's transitivity theorem for JB$^*$-triples \cite[Theorem 3.3]{BurFerGarMarPe} or by its version for JB$^*$-algebras \cite[Proposition 2.3]{Ham99}, there exist mutually orthogonal norm-one positive elements $a_1,\ldots,a_k$ in $\AAA$ such that $a_j = P_2(p_j) (a_j) + P_0(p_j) (a_j)$ for all $j=1,\ldots, k$. By orthogonality, it is easy to check that $h = \sum_{j=1}^{k} \theta_j a_j$ is a hermitian element in $\AAA$ with $h = P_2(p_j) (h) + P_0(p_j) (h)$ for all $j=1,\ldots, k$. Set $p_s := \sum_{j=1}^k p_j\in \AAA^{**}.$  By orthogonality and the properties of the $a_j$'s we have $$h = \sum_{j=1}^k \theta_j p_j + P_0(p_s) (h), \hbox{ and } \exp(i h) = \sum_{j=1}^k \exp(i \theta_j) p_j + P_0(p_s) (\exp(ih)),$$ where the summands in these sums are mutually orthogonal. Clearly, $e\leq u \leq \exp(i h)$ by construction. Clearly the unitary $\exp(ih)$ lies in the principal connected component of $\mathcal{U} (\AAA)$ (cf. \eqref{eq algebraic charact principal component unitaries JBstar}). 
\end{proof}

We are now in a position to obtain an strengthened version of Theorem~\ref{Russo_Dye_theorem_for_F_p}. 

\begin{theorem}\label{Russo_Dye_theorem_for_maximal proper faces} Let $F$ be a maximal proper faces of the closed unit ball of a unital JB$^*$-algebra $\AAA$. Then $F=\overline{\mathrm{co}}\left(F\cap{\mathcal U}^0(\AAA)\right)$.
\end{theorem}

\begin{proof} Since $F$ is a maximal proper face of $\mathcal{B}_{\AAA}$, there exits a minimal tripotent $e\in \AAA^{***}$ satisfying $F = F_e$ (cf. Theorem~\ref{t facil structure of the ball of a JB$^*$-triple} or \cite[Corollary 3.5 and the preceding comments]{BuFerMarPe}). By Theorem~\ref{t boundedness  principle for minimal tripotents in the second dual by unitaries in the principal component} there exists a unitary $u$ in the principal connected component of $\mathcal{U} (\AAA)$ such that $e\leq u$. Having in mind that $\mathcal{U}^0 (\AAA(u))  =\mathcal{U}^0 (\AAA)$ and the fact that $e$ is a minimal projection in the bidual of the $u$-isotope $\AAA(u)$, the desired conclusion is a direct consequence of Theorem~\ref{Russo_Dye_theorem_for_F_p}.
\end{proof}

The last result of this section is an strengthened version of Proposition~\ref{Prop_20_from_Mori-Ozawa}.

\begin{theorem}\label{theorem prop_20_from_Mori-Ozawa for maximal faces associated with minimal tripotents}
	Let $\AAA$ be a unital JB$^*$-algebra, $Y$ a Banach space, $\Delta:S_\AAA\to S_Y$ a surjective isometry, and $F_e\subset S_{\AAA}$ a maximal norm-closed proper face associated with a minimal tripotent $e\in \AAA^\dd$. Then there exist a directed set $\Lambda$ and a net $(\Theta_\lambda)_\lambda$ of affine contractions $\Theta_\lambda: F_e\to K_\lambda,$ where each $K_\lambda$ is a closed convex subset of $F_e$, such that $\Theta_\lambda\to\mathrm{id}_{F_e}$ in the point-norm topology {\rm(}i.e. in the strong operator topology{\rm)}, each $\Delta(K_\lambda)$ is convex, and each $K_\lambda$ is affinely isometrically isomorphic to the closed ball of a real JB$^*$-algebra $\mathcal{R}_\lambda$ such that $\B_{\mathcal{R}_\lambda}$ satisfies the strong Mankiewicz property {\rm(}and therefore $\Delta\zuz_{F_e}$ is affine{\rm )}. In particular $\Delta$ is affine on every norm closed proper face of $\mathcal{B}_{\AAA}$.   
\end{theorem}

\begin{proof} Let $F_e$ be a maximal proper face of $\mathcal{B}_{\AAA}$ associated with a minimal tripotent $e\in \AAA^{**}$. By a new application of Theorem~\ref{t boundedness  principle for minimal tripotents in the second dual by unitaries in the principal component} there exists a unitary $u$ in $\AAA$ such that $e\leq u$. This implies that $e$ is a minimal projection in the bidual of the $u$-isotope $\AAA (u)$. Therefore the desired conclusion follows from Proposition~\ref{Prop_20_from_Mori-Ozawa}.
\end{proof}

\section{Applications to Tingley's problem and the Mazur--Ulam property}\label{sec: proof of the main theorem}

The algebraic--geometric results developed in the previous sections are interesting tools by themselves with many potential applications. For the purposes of this paper we shall apply them to provide a complete positive solution to Tingley's problem in the case in which one of the spaces is a unital JB$^*$-algebra. As we shall see later, we can restrict our discussion to the case of unital JB$^*$-algebras of infinite rank, so the technical auxiliary results are in this context.\smallskip

We recall that, according to \cite[Definition 4.8]{BeCuFerPe2018}, a JB$^*$-triple $E$ satisfies property {\rm(}$\mathcal{P}${\rm)} if for each minimal tripotent $e$ in $E^{**}$ and each complete tripotent $u$ in $E$ {\rm(}equivalently, $u\in \partial_e(\mathcal{B}_{E})${\rm)}, there exists another minimal tripotent $w$ in $E^{**}$ satisfying $w\perp e$ and $u = w + P_0(w) (u)$. It is shown in the proof of \cite[Theorem 4.14]{BeCuFerPe2018} that every JBW$^*$-triple with rank bigger than or equal to three satisfies property {\rm(}$\mathcal{P}${\rm)}. Proposition 4.13 in \cite{BeCuFerPe2018} shows that every Cartan factor of rank bigger than or equal to three enjoys this property, however the proof is literally valid to get the following lemma whose proof is omitted. 

\begin{lemma}\label{l every Cartan factor satisfies Pprime}{\rm(\cite[Proposition 4.13]{BeCuFerPe2018})} Let $C$ be a Cartan factor of rank bigger than or equal to three. Then for each minimal tripotent $v$ in $C$ and each complete tripotent $u$ in $C$ there exists another minimal tripotent $e$ in $C$ satisfying $v\perp e$ and $u = e + P_0(v) (u)$.
\end{lemma}

We can now mimic the arguments in the proof of \cite[Theorem 4.14]{BeCuFerPe2018} to derive an strengthened version of the conclusion.

\begin{lemma}\label{L principal_unitary_lies_in_a_perpendicular_face} Let $E$ be a JB$^*$-triple of infinite rank, $u$ a complete tripotent in $E$ and $v\in E^\dd$ a minimal tripotent. Then there is a minimal tripotent $e$ in $E^\dd$ such that $e\perp v$ and $u\in F_e$.
\end{lemma}

\begin{proof} Clearly, the atomic part, $\mathcal{A},$ of $E^{**}$ is a JBW$^*$-triple of infinite rank since, as we commented at the beginning of the previous section, $E$ embeds isometrically into this atomic part. \smallskip
	 
We can write the atomic part of $E^{**}$ in the form $\mathcal{A} = \bigoplus^{\ell_{\infty}}_{j\in\Lambda} C_j$, where each $C_j$ is Cartan factor  \cite[Proposition 2]{FriRu86}.	If $\iota: E\hookrightarrow E^{**}$ and $\pi_{\mathcal{A}}: E^{**}\to \mathcal{A}$ stand for the canonical embedding of $E$ into its bidual and the projection of the latter onto its atomic part, respectively, the mapping $ \pi_{\mathcal{A}}\circ \iota : E\to \mathcal{A}$ is an isometric triple embedding. 
 The element $u$ also is a complete tripotent in $E^{**},$ just as $\pi_{\mathcal{A}} (u)$ is  a complete tripotent in $\mathcal{A}$. Assuming that $\sharp \Lambda \geq 2,$ by minimality there exists $j_1\neq j_2$ in $\Lambda$ such that $v\in C_{j_1}$ and the $j_2$-component of $ \pi_{\mathcal{A}}(u)$ is a non-zero complete tripotent in $C_{j_2}$. Since every non-zero tripotent in $C_{j_2}$ is an $\ell_{\infty}$-sum of a family of mutually orthogonal minimal tripotents, we can find a minimal tripotent  $e$ in $C_{j_2}$ (and also in $E^{**}$)  with $e\leq \pi_{\mathcal{A}} (u)\le u$ and $e\perp v$. Clearly, $u\in F_e$. \smallskip

If $\sharp \Lambda =1,$ the atomic part of $E^{**}$ reduces to a single Cartan factor of infinite rank. Lemma \ref{l every Cartan factor satisfies Pprime} applied to $\pi_{\mathcal{A}} (u)$ proves the desired statement. 
\end{proof}

In our next result we combine a non-commutative version of Urysohn's \cite{FerPeUrysohn} with a recent property of norm closed proper faces. We recall first that the orthogonal complement of an element $b$ in a JB$^*$-triple $E$ is defined as $\{b\}^{\perp}_{_E} \df \{x\in E : x\perp b\}$. It is known that $\{b\}^{\perp}_{_E}$ is an inner ideal of $E$ which coincides with $E^{**}_0 (r(b)) \cap E$ (see, for example, \cite[Lemma 3.11]{BurGarPe11} and \eqref{eq reformulations of orthogonality}). 

\begin{lemma}\label{L minimal_tripotent_is_in_weak*_closure}
Let $E$ be a JB$^*$-triple and let $e$ be a minimal tripotent in $E^\dd$. Let $b$ be an element in $\AAA$ satisfying $b\perp e$ and 
$F_e \cap \{b\}_{_E}^{\perp}\neq \emptyset.$ Then $e\in\overline{F_e \cap \{b\}_{_E}^{\perp}}^{w\ast}$.
\end{lemma}

\begin{proof} If $b=0$ the set $F_e \cap \{b\}_{_E}^{\perp}$ coincides with $F_e = \left(e + \mathcal{B}_{E^{**}_0(e)} \right)\cap E$, and we know from \cite[Theorem 3.6]{BeCuFerPe2018} that $$\overline{F_e}^{w^*} = F_e^{^{E^{**}}} =  e + \mathcal{B}_{E^{**}_0(e)},$$  which clearly contains $e$. \smallskip
	
We assume next that $b\neq 0$. By hypothesis, we can take $a\in F_e \cap \{b\}_{_E}^{\perp}$. It is known that the inner ideal of $E$ generated by $a,$ $E(a)$, coincides with the norm closure of $Q(a) (E)$, its weak$^*$-closure in $E^{**}$, which clearly identifies with $E(a)^{**}$, is precisely $E_2^{**} (r(a))$, and $E(a)$ is a JB$^*$-subalgebra of $(E^{**}_2 (r(a)), \circ_{r(a)}, *_{r(a)})$ (cf. \cite[Proposition 2.1 and comments prior to it]{BunChuZal2000}). By Peirce arithmetic, or by \eqref{eq reformulations of orthogonality}, $E(a)\subseteq \{b\}^{\perp}_{_E}.$\smallskip

Since $a\in F_e$, we easily deduce that $r(a)\geq e$ in $E^{**}$, and hence $e\in E^{**}_2 (r(a)) = E(a)^{**}$. One of the consequences of the non-commutative version of Urysohn's lemma proved in \cite[Corollary 4.2]{FerPeUrysohn} asserts that $e$ is a compact tripotent in $E(a)^{**}$. To simplify the notation, let us write $I = E(a)$. According to this notation, let $F_e^{^I}$ and $F_e^{^{I^{**}}}$ denote the norm closed proper face of $\mathcal{B}_I$ and the weak$^*$ closed proper face of $\mathcal{B}_{I^{**}}$ associated with the compact tripotent $e$, respectively. Clearly, $F_e^{^I} = F_e\cap I$, and it follows from \cite[Theorem 3.6]{BeCuFerPe2018} that $$\overline{F_e\cap {I}}^{w^*} = \overline{F_e^{^I}}^{w^*} = F_e^{^{I^{**}}} =  e + \mathcal{B}_{I^{**}_0(e)}\ni e,$$ which proves the desired statement since $F_e \cap I \subseteq F_e\cap \{b\}^{\perp}_{_E}$. 
\end{proof}

We recall some useful results. For each Banach space $X$, the maximal convex subsets of $S_{X}$ are precisely the maximal proper norm closed faces of $\mathcal{B}_{X}$ (cf. \cite[Lemma 3.3]{Tan2016} or \cite[Lemma 3.2]{Tan2017b}). Furthermore, for each maximal proper norm closed face $F$ of $\mathcal{B}_X,$ there exists an extreme point $\varphi$ of the closed unit ball $\mathcal{B}_{X^*}$ such that $F = \varphi^{-1}\{1\}\cap  S_{X}$ (cf. \cite[Lemma 3.3]{Tan2016}). The set of all extreme points $\varphi$ of  $\mathcal{B}_{X^*}$ for which $\varphi^{-1}\{1\}\cap S_{X}$ is a maximal convex subset of $S_{X}$ will be denoted by $\mathcal{Q}_{X}.$ Since, by Zorn's lemma, each $x\in S_{X}$ is contained in a maximal convex subset of $S_{X}$, the set $\mathcal{Q}_{X}$ is norming on $X$, actually the norm of each element in $X$ is attained at an element in $\mathcal{Q}_{X}$.

\begin{proposition}\label{P phi_circ_Delta_=_varphi} 
Let $\AAA$ be a unital JB$^*$-algebra of infinite rank, $Y$ a real Banach space, $\Delta: S_{\AAA}\to S_{Y}$ a surjective isometry, and $\phi\in \mathcal{Q}_{Y}$. Then there exists $\varphi\in\partial_e(\B_{\AAA^\ast})$ such that $$\phi\circ\Delta=\Re\varphi\zuz_{S_{\AAA}}.$$
\end{proposition}

\begin{proof} Let $\mathbf{F}\df\phi^{-1}(\{1\})$. By hypotheses, $\mathbf{F}$ is a maximal proper norm closed face of $\mathcal{B}_Y$ (equivalently, a maximal convex set of $S_{Y}$). By \cite[Lemma 5.1]{ChenDong2011} (or \cite[Lemma 3.5]{Tan2014}), $\Delta^{-1} (\mathbf{F})$ is a maximal convex subset of $S_{\AAA}$, or what is equivalent, a maximal proper norm closed face of $\mathcal{B}_{\AAA}$. Therefore, there exists a minimal tripotent $v\in \AAA^{**}$ such that $F_v = \Delta^{-1} (\mathbf{F})$ (cf. Theorem \ref{t facil structure of the ball of a JB$^*$-triple}). \smallskip

Therefore, there is a unique pure atom $\varphi_v\in\partial_e(\mathcal{B}_{\AAA^\ast})$ such that $F_v = \varphi_v^{-1} (\{1\}),$ and hence $\varphi_v\zuz_{F_v}\equiv 1$. We claim that
$$\hbox{ $\phi\circ\Delta=\Re\hbox{e} \varphi_v$ on the whole $S_\AAA$.	}$$
Having in mind Theorems \ref{Russo_Dye_theorem_for_maximal proper faces} and \ref{theorem prop_20_from_Mori-Ozawa for maximal faces associated with minimal tripotents}, it suffices to prove that $\phi\circ\Delta(u)=\varphi_p(u)$ for all $u\in\mathcal{U}^0(\AAA)$. Fix an arbitrary unitary  $u$ in $\AAA$.\smallskip

By Lemma \ref{L principal_unitary_lies_in_a_perpendicular_face} there is a minimal tripotent $e\in\AAA^\dd$ such that $e\perp v$ and $u\in F_e$. 	Since, by Theorem \ref{theorem prop_20_from_Mori-Ozawa for maximal faces associated with minimal tripotents}, $\Delta$ is affine on $F_e$, there is $\gamma_e\in\R$ and $\widetilde{\varphi}_e\in \mathcal{B}_{\AAA^\ast}$ such that $\phi\circ\Delta(y)=\gamma_e+\Re\hbox{e} \widetilde{\varphi}_e (y)$ for any $y\in F_e$. \smallskip

Using Kadison's transitivity theorem for JB$^*$-triples \cite[Theorem 3.3]{BuFerMarPe} we can find $a\in F_v$ and $b\in F_e$ such that $a\perp b$. So, the element $a$ lies in the set $F_v \cap\{b\}^{\perp}_\AAA,$ which is thus non-empty. Fix any other $x\in F_v \cap\{b\}^{\perp}_\AAA$. Then we have $b\pm x\in F_e\cap\left(\pm F_v\right)$, and hence, by the above properties, we arrive at $$\pm1 =\phi\circ\Delta(b\pm x)=\gamma_e+  \Re\hbox{e} \widetilde{\varphi}_e  (b\pm x).$$ This implies that  $\gamma_e+ \Re\hbox{e} \widetilde{\varphi}_e (b)=0$ and $\Re\hbox{e} \widetilde{\varphi}_e (x)=1$. Therefore we derived that $\Re\hbox{e} \widetilde{\varphi}_e (x)=1$ for every $x\in F_v \cap\{b\}^{\perp}_\AAA$.\smallskip
	
Now, Lemma \ref{L minimal_tripotent_is_in_weak*_closure} asserts that $v\in\overline{F_v \cap\{b\}^{\perp}_\AAA}^{w*}$, so we must have $ \Re\hbox{e} \widetilde{\varphi}_e (v)=1$, and $\|\Re\hbox{e} \widetilde{\varphi}_e \|\le 1$ implies that $\Re\hbox{e} \widetilde{\varphi}_e =\Re\hbox{e} \varphi_v$ (see the discussion on page \pageref{eq extreme points and complete tripotents}). Finally, since $b\in F_e$ and $b\perp v$ we also have  $0=\gamma_e+ \Re\hbox{e} \varphi_v (b)=\gamma_e$. Therefore $\phi\circ\Delta= \Re\hbox{e} \varphi_v$ on $F_e$, in particular $\phi\circ\Delta(u)=  \Re\hbox{e} \varphi_v (u)$.
\end{proof}

We are now in a position to prove that every unital JB$^*$-algebra satisfies the Mazur--Ulam property. 

\begin{theorem}\label{t unital JB$^*$-algebras satisfy the MUP} Every unital JB$^*$-algebra $\AAA$ satisfies the Mazur--Ulam property, that is, for each Banach space $Y$, every surjective isometry $\Delta: S_{\AAA}\to S_{Y}$ admits an extension to a surjective real linear isometry from $\AAA$ onto $Y$. 
\end{theorem}

\begin{proof} If $\AAA$ has finite rank it must be reflexive (see, for example, \cite[Proposition 4.5]{BunChu92} and \cite[Theorem 6]{ChuIo90}). In particular, $\AAA$ is a JBW$^*$-triple, and hence the desired result follows from \cite[Corollary  1.2]{KalPe2019} even in a more general setting. We can therefore assume that $\AAA$ (and hence $\AAA^{**}$) has infinite rank.\smallskip  
	
Having in mind \cite[Lemma 2.1]{FangWang06}, the desired conclusion holds if we show that $\|\Delta(a)-\lambda\Delta(b)\|\le\|a-\lambda b\|$ for all $a,b\in S_\AAA$ and $\lambda>0$. Fix $a,b,$ and $\lambda$ under these conditions. Since every element in $Y$ attains its norm at an element in  $\mathcal{Q}_{Y}$, there exists $\phi\in \mathcal{Q}_{Y}$ such that $\phi(\Delta(a)-\lambda\Delta(b))=\|\Delta(a)-\lambda\Delta(b)\|$. By Proposition \ref{P phi_circ_Delta_=_varphi} we can find a pure atom $\varphi\in \partial_e \left( \mathcal{B}_{\AAA^\ast}\right)$ such that $\phi\circ\Delta=\Re\hbox{e}\varphi\zuz_{S_{\AAA}}$. So \begin{align*}
		\|\Delta(a)-\lambda\Delta(b)\|&=\phi(\Delta(a)-\lambda\Delta(b))=\phi(\Delta(a))-\lambda \phi(\Delta(b))\\
		&=\Re\hbox{e}\varphi(a)-\lambda\Re\hbox{e}\phi(b)=\Re\hbox{e}\varphi(a-\lambda b)\le \|a-\lambda b\|,
	\end{align*} which concludes the proof.	
\end{proof}

\medskip\medskip

\textbf{Acknowledgements} A.M. Peralta partially supported by grant PID2021~-122126NB--C31 funded by MCIN/AEI/10.13039/501100011033 and by ``ERDF A way of making Europe'', Junta de Andaluc\'{\i}a grants FQM375 and PY20$\underline{\ }$00255, and by the IMAG--Mar{\'i}a de Maeztu grant CEX2020-001105-M/AEI/10.13039/ 501100011033.  R. \v{S}varc is supported by the Charles University, project GAUK no. 268521. Part of this research was conducted during a research visit of the second author to the Department of Mathematical Analysis and the Math Institute of the University of Granada; he would like to express his gratitude for the hospitality.

\end{document}